\documentclass[12pt]{article}

\usepackage[top=1in, bottom=1in, left=1.5in, right=1.5in]{geometry}

\usepackage[dvipsnames, table, xcdraw]{xcolor}

\usepackage[toc,page]{appendix}

\usepackage[utf8]{inputenc}        
\usepackage{sans} 
\usepackage{pifont} 

\usepackage{graphicx}              
\usepackage{caption,subcaption}    
\usepackage{hyperref}              
\hypersetup{                       
	colorlinks=true,               
	linkcolor=blue,                
	filecolor=magenta,             
	urlcolor=cyan,                 
	citecolor=ForestGreen          
}
\graphicspath{{Figures/PDFs}}

\usepackage{amsmath}               
\usepackage{amssymb}               
\usepackage{amsfonts}              
\usepackage{amsthm}                
\usepackage{mathtools}
\usepackage{array}
\usepackage{multirow}
\usepackage{booktabs}

\usepackage{enumitem}  
\usepackage[breakable,skins]{tcolorbox}

\theoremstyle{style-one}
\newtheorem{theorem}{Theorem}[section]  		

\newtheorem*{definition*}{Definition}  
\newtheorem*{remark}{Remark}

\newtheorem*{problem*}{Problem}  

\newcommand{\uP}{u}
\newcommand{\sP}{\mathbf{s}}
\newcommand{\eP}{\mathbf{e}}
\newcommand{\lP}{\lambda}
\newcommand{\mP}{\boldsymbol{\mu}}

\newcommand{\stP}{\mathbf{\tilde{s}}}
\newcommand{\etP}{\mathbf{\tilde{e}}}

\newcommand{\CG}{\text{CG}}

\newcommand{\DG}{\text{DG}}

\usepackage{textcomp, upquote}

\usepackage{listings}

\definecolor{codegreen}{rgb}{0,0.6,0}
\definecolor{codegray}{rgb}{0.5,0.5,0.5}
\definecolor{codepurple}{rgb}{0.58,0,0.82}
\definecolor{backcolour}{rgb}{0.95,0.95,0.92}

\lstdefinestyle{mystyle}{
    backgroundcolor=\color{backcolour},
    commentstyle=\color{codegreen},
    keywordstyle=\color{magenta},
    numberstyle=\tiny\color{codegray},
    stringstyle=\color{codepurple},
    basicstyle=\ttfamily\scriptsize,
    breakatwhitespace=false,
    breaklines=true,
    captionpos=b,
    keepspaces=true,
    numbers=none,
    numbersep=5pt,
    showspaces=false,
    showstringspaces=false,
    showtabs=false,
    tabsize=2,
    upquote=true
}

\lstnewenvironment{python}[1][]
{
\pythonstyle
\lstset{#1}
}
{}

\newcommand\pythoninline[1]{{\pythonstyle\lstinline!#1!}}

\usepackage{authblk}

\title{Numerical approximation of a PDE--constrained Optimization problem that appears in  Data-Driven Computational Mechanics}

\date{}

\author[1]{Pedro B. Bazon}

\author[2*]{Cristian G. Gebhardt}

\author[1]{Gustavo C. Buscaglia}

\author[1]{Roberto F. Ausas}

\affil[1]{Instituto de Ci\^encias Matem\'aticas e de Computa\c{c}\~ao, Universidade de S\~ao Paulo\\ Av. Trab. S\~ao Carlense 400, 13566-590, São Carlos-SP, Brazil}

\affil[2]{University of Bergen, Geophysical Institute and Bergen Offshore Wind Centre\\ All\'{e}gaten 70, 5007 Bergen, Norway.}

\begin{document}

\maketitle

\let\thefootnote\relax\footnotetext{\textsuperscript{*}Corresponding author email: cristian.gebhardt@uib.no} 

\begin{abstract}

We investigate an optimization problem that arises when working within the paradigm of Data-Driven Computational Mechanics. In the context of the diffusion-reaction problem, such an optimization problem seeks for the continuous primal fields (gradient and flux) that are closest to some predefined discrete fields taken from a material data set. The optimization is performed over primal fields that satisfy the physical conservation law and the geometrical compatibility. We consider a reaction term in the conservation law, which has the effect of coupling all the optimality conditions. We first establish the well-posedness in the continuous setting. Then, we propose stable finite element discretizations that consistently approximate the continuous formulation, preserving its saddle-point structure and allowing for equal-order interpolation of all fields. Finally, we demonstrate the effectiveness of the proposed methods through a set of numerical examples.

\end{abstract}



\section{Introduction}

\textit{Data-Driven Computational Mechanics} (DDCM) represents the latest paradigm in continuum mechanics and can be employed in a wide range of applications. 
For the diffusion-reaction problem, for example, DDCM seeks to find the gradient of a concentration scalar field (a conserved quantity) and the flux vector field that are closest to a field extension of empirically determined gradient-flux pairs.
Such a problem can be stated as an optimization problem in which discrete and continuous fields are to be simultaneously sought such that they are as close as possible while the continuous ones satisfy the physical conservation law and the compatibility geometrical condition. 

In this context, a basic question that arises is about the existence of optimal continuous fields {\em when the discrete fields are fixed and given}. This ``optimization subproblem'' is actually a continuous optimization problem that can be treated computationally and analyzed using standard tools from the classical theory of linear \textit{partial differential equations} (PDEs). This is the main target of the present contribution. To place the proposed ideas into an appropriate setting, we provide below some of the most important contributions to present date, highlight the difference of the present work with respect to existing literature and state the novelty from the methodological point of view.

DDCM is a very active field of research whose first contributions can be traced back to less than a decade ago and has mostly been developed within solid/structural mechanics. For instance, this principle is applied to the static analysis of multi-bar structures \cite{Kirchdoerfer2016}. Based on similar ideas, a data-driven magnetostatic solver is proposed in \cite{Galetzka2020}. In the following, this initial approach is modified to improve robustness due to the presence of outliers in the constitutive data \cite{Kirchdoerfer2017}. Finally, this is extended to the dynamic analysis of structures that exhibit configuration-independent mass matrices \cite{Kirchdoerfer2018}. For the sake of completeness, we ought to mention that there are some works addressing the existence of solutions for data-driven elasticity \cite{Conti2018,Conti2020}. There is also a hybrid approach that relies on a smooth constitutive manifold reconstructed from the data at an offline pre-processing step. This principle is applied to time-independent problems \cite{Gebhardt:2020a}. Then, the approach is also modified to handle time-dependent problems \cite{Gebhardt:2020b}. Lastly, this is extended to deal with \textit{mathematical programs with complementarity constraints} (MPCC) such as those arising in contact mechanics  \cite{Gebhardt2024a}. Recently in  \cite{Gebhardt2024b}, the structural properties of the DDCM problem discretized by means of the \textit{finite element method} (FEM) is investigated and its global solvability is established. Remarkably, there are neither antecedents of DDCM in the context of the diffusion-reaction problem nor direct connections to PDE-regularized methods to analyze empirical data \cite{Ponti2021} or to inhomogeneities identification from boundary measurements \cite{Ammari2004}.

All the works introduced above (that deal with the numerical solution of the discrete-continuous optimization problem in DDCM) require along the optimization process to compute the solution of a type of continuous ``optimization subproblem" that is a particular case of the optimization problem addressed in the present work. This is because those works adopt an \textit{alternating direction method} (ADM). This work can be considered then as a first step towards the application of DDCM in the context of the diffusion-reaction problem, which to the best of our knowledge has not been realized so far. In previous works, once the assignation of the data pairs has been done, the optimization subproblem decomposes into two sets of optimality conditions that are totally decoupled and can be solved independently. In the present work, we consider a reaction term present in the physical conservation law. This automatically couples all the optimality conditions and thus leads to a novel mathematical problem, for which we provide a first well-posedness analysis. Moreover, our analysis considers two and three spatial dimensions in contrast to contributions that focus on problems defined on closed intervals of the real line \cite{Gebhardt:2020a, Gebhardt:2020b, Gebhardt2024a, Gebhardt2024b}. In addition, we introduce new finite element formulations, establish their well-posedness for the chosen finite element spaces and show their convergence to the exact solution. The desired stability properties are achieved by directly approximating the saddle-point problem, approach that exhibits some connection to other well-established approaches \cite{Masud2002, Bochev2006, Badia2010, Burman2023}. 

The remaining of this article is organized as follows. In Section 2, we describe the optimization problem in the reaction-diffusion context. We prove that this optimization problem, in its continuous setting, is indeed well posed. We also propose and analyze finite element approximations of the continuous problem. In Section 3, we introduce the numerical implementation of the proposed methods, which is realized by coding within the Firedrake finite element platform. 
This development is followed as well by a series of examples, of increasing difficulty, that illustrate the remarkable properties of the setting proposed. Finally, in Section 4, we conclude with remarks, potential limitations and possible future work.

\section{The model problem}

We consider a reaction-diffusion problem governed by the conservation law
\begin{equation}
	\nabla \cdot \sP + \zeta\,\uP = q
	\label{eq:conservation}
\end{equation}
posed on a regular domain $\Omega \subset \mathbb{R}^{d_{\Omega}}$. In eq.~(\ref{eq:conservation}) the symbol $\nabla$ represents the vector of partial derivatives, so that $\nabla \cdot \sP$ is the divergence of the vector field $\sP$, which represents the flux. Also, $\uP$ is the concentration of the conserved quantity, $\zeta \in L^{\infty}(\Omega)$ is the (possibly position-dependent) reaction coefficient and $q\in L^2(\Omega)$ is the source. In thermal applications $\uP$ is the temperature, $\sP$ is the heat flux and $q-\zeta \uP$ is a heat source that depends linearly on the temperature. In weakly compressible flow in porous media $\uP$ is the pressure and $\sP$ is the fluid flux.

The boundary of $\Omega$, denoted by $\partial \Omega$, is divided into a Dirichlet part $\partial \Omega_D$, where $\uP=\uP_D$ is imposed, and a Neumann part $\partial\Omega_N$, in which $\sP \cdot \mathbf{\check{n}}=H$ is imposed. Hereafter, we assume for simplicity homogeneous boundary conditions throughout, i.e., $\uP_D = 0$ and $H = 0$.

Let us also define the field $\eP$ as the gradient of $\uP$,
\begin{equation}
	\eP = \nabla \uP
\label{eq:compatibility}
\end{equation}

This equation can be viewed as a compatibility condition. We call a vector field $\eP$ compatible if there exists a scalar field $\uP$ such that $\eP=\nabla \uP$.

Now, suppose that two vector fields $\stP$ and $\etP$ are given in $\Omega$. These fields will both be assumed in $\left[L^2(\Omega)\right]^{d_{\Omega}}$ and thus possibly discontinuous. Since they originate from a measured dataset, it is generally the case that $\stP$ does not satisfy the conservation law, and $\etP$ does not satisfy the compatibility condition. Moreover, they are concomitant, in the sense that each sample $\etP$ corresponds to the same spatial location as $\stP$ -- they form coupled observations of gradient and flux at the same point in the domain --.

The goal, therefore, is to determine the best approximations to $\stP$ and $\etP$ within the space of feasible fields, that is, those that are both conservative and compatible. The scalar field $\uP$ emerges naturally as a by-product of this reconstruction.

\subsection{Well-posedness of the Exact Problem}

We begin by introducing the function spaces involved in the formulation:

\begin{equation}
	\left.
	\begin{array}{rcl}
		U & = & \left\{ \bar{u} \in H^1(\Omega),\ \bar{u} = 0 \text{ on } \partial\Omega_D \right\} \\
		E & = & \left[ L^2(\Omega) \right]^{d_\Omega} \\
		S & = & \left[ L^2(\Omega) \right]^{d_\Omega} \\
	\end{array}
	\right \rbrace  \quad X := U \times E \times S
\end{equation}
We now pose the problem as a minimization over the feasible set $ Z \subset X $, consisting of all tuples $ (\uP, \eP, \sP) \in X $ that satisfy the governing equations and boundary conditions.
\begin{tcolorbox}[colback=black!10!white, halign = justify, boxrule=0pt]
\begin{problem*}[P1]
	Find $x= \left (\uP,\eP,\sP \right) \in Z$ -- the feasible space -- that minimizes the distance to the data in the feasible space, i.e.,
	\begin{equation}
		J(x) \leq J(x+\delta x)~,\qquad \forall \,\delta x\,\in\,\delta Z~
	\end{equation}
	where 
	\begin{equation}
		J(\bar{\uP},\bar{\eP},\bar{s})=\frac12 \|\bar{\eP}-\etP\|^2 + \frac{\kappa}{2} \|\bar{\sP}-\stP\|^2~,
	\end{equation}
	and the feasible space is
	\begin{equation}
		Z=\left \{ (\bar{\uP},\bar{\eP},\bar{\sP})\in X~|~\nabla \cdot \bar{\sP} + \zeta \,\bar{\uP} = q~,~\bar{\eP}=\nabla \bar{\uP},\,\left.~\bar{\sP}\cdot\mathbf{\check{n}}\right |_{\partial \Omega_N}=0 \right \}
	\end{equation}
	with the variation given by
	\begin{equation}
		\delta Z=\left \{ (\delta{\uP},\delta{\eP},\delta{\sP})\in X~|~\nabla \cdot \delta{\sP} + \zeta \delta{\uP} = 0~,~ \delta{\eP}=\nabla \delta{\uP} ,\,\left. ~\delta{\sP}\cdot\mathbf{\check{n}}\right |_{\partial \Omega_N}=0\right \}~.
	\end{equation}
\end{problem*}
\end{tcolorbox}

\begin{remark}
	Above, $\|\cdot\|$ is the norm of $\left[L^2(\Omega)\right]^{d_{\Omega}}$. The corresponding scalar product is denoted by $\left ( \cdot , \cdot \right )$. The same notation will be used for $L^2(\Omega)$.
\end{remark}
Let us now relax the constraints through Lagrange multipliers $\lP$ and $\mP$, which belong to the following function spaces:
\begin{equation*}
	L=\{g\in H^1(\Omega)~|~g=0 \mbox{ on } \partial \Omega_D\},\qquad
	\mbox{ and }\qquad M=\left[L^2(\Omega)\right]^{d_{\Omega}},  
\end{equation*}
respectively. We also use the notation
\begin{equation}
	y=\left (\lP,\mP\right )\,\in \, Y:=\,L\times M~.
\end{equation}
Then, the resulting equivalent problem can be written as: 
\begin{tcolorbox}[colback=black!10!white, halign = justify, boxrule=0pt]
	\begin{problem*}[P2]
		Find $(x,y)=\left( \uP, \eP, \sP, \lP, \mP \right) \in X\times Y$ such that 
		\begin{equation}
			(x,y)= \arg \left \{ \inf_{\bar{x}\in X} \sup_{\bar{y}\in Y} 
			\mathcal{L}(\bar{x},\bar{y})
			\right \}
		\end{equation}
		where
		\begin{equation}
			\begin{aligned}
				\mathcal{L}(\bar{x},\bar{y}) &~=~J(\bar{\uP},\bar{\eP},\bar{\sP})+ \nonumber \\
				&~-~\left ( \bar{\sP},\nabla \bar{\lP} \right ) + \left (\zeta \bar{\uP} - q , \bar{\lP} \right )  
				+ \left ( \bar{\eP}-\nabla \bar{\uP}, \bar{\mP} \right ) ~.
			\end{aligned} \label{eq:lagrangian}
		\end{equation}
	\end{problem*}
\end{tcolorbox}
The weak form of the optimality conditions reads:
\begin{equation}\label{eq:variational_form}
	\begin{array}{rcll}
		\left ( \lP, \zeta \delta \uP \right ) - \left ( \mP, \nabla \delta \uP \right ) & = & 0 &\qquad \forall \,\delta \uP \in U \\
		\left ( \eP, \delta \eP \right ) + \left ( \mP, \delta \eP \right ) &=& \left ( \etP, \delta \eP \right )  &\qquad \forall \,\delta \eP \in E\\
		\kappa \left ( \sP, \delta \sP \right ) - \left ( \nabla \lP, \delta \sP \right )  & = & \kappa \left ( \stP,\delta \sP \right ) &\qquad \forall\, \delta \sP \in S \\
		\left ( \zeta \,\uP, \delta \lP \right )-\left ( \sP, \nabla \delta \lP \right )  &=& \left ( q, \delta \lP \right ) 	&\qquad \forall\,\delta \lP \in L\\
		-\left (\nabla \uP,\delta \mP \right ) + \left ( \eP, \delta \mP \right ) & = & 0 &\qquad \forall \, \delta \mP \in M 
	\end{array}
\end{equation}

The sum of the left-hand sides (with the last two equations multiplied by $-1$) defines the bilinear form of the problem in the product space $X \times Y$, given by:

\begin{equation}
	B((x,y),(\delta x,\delta y))=a(x,\delta x)+b(\delta x,y)-b(x,\delta y)
\end{equation}
where
\begin{equation}\label{eq:bform}
	\begin{array}{rcl}
		a(x,\delta x)&=& \left (\eP,\delta \eP \right )+\kappa \left ( \sP,\delta \sP \right ) \\
		b(x,y)& = & \left ( \mP,\eP-\nabla \uP \right )-\left (\nabla \lP,\sP\right ) + \left ( \lP, \zeta \uP \right )
	\end{array}
\end{equation}

The optimality conditions in strong form, therefore, are given by:
\begin{eqnarray}
	\zeta \lP  + \nabla \cdot \mP & = & 0 			\label{eq:dif1}\\
	 \eP + \mP & = & \etP 							\label{eq:dif2}\\
	\kappa\, \sP - \nabla \lP & = & \kappa\,\stP 	\label{eq:dif3}\\
	 \zeta\,\uP+\nabla \cdot \sP  & = & q 			\label{eq:dif4}\\
	-\nabla \uP + \eP &=& 0 						\label{eq:dif5}
\end{eqnarray}
together with
\begin{eqnarray}
 \sP \cdot \mathbf{\check{n}} & = & 0~,\qquad \mbox{ on } \partial \Omega_N \\
 \mP \cdot \mathbf{\check{n}} & = & 0~,\qquad \mbox{ on } \partial \Omega_N ~.
\end{eqnarray}

\begin{remark}
	Notice that when the reaction coefficient $\zeta$ is identically zero the five variational equations given in (\ref{eq:variational_form}) -- or, equivalently, their differential versions, i.e., equations (\ref{eq:dif1}) to (\ref{eq:dif5}) -- are in fact two independent sets of equations: (a) The first, second and fifth equations, with unknowns $\uP$, $\eP$ and $\mP$, and (b) the third and fourth equations, with unknowns $\sP$ and $\lP$. In such a case, $\sP$ is the best approximation to $\stP$ among conservative flux fields, and $\eP=\nabla \uP$ is the best approximation to $\etP$ among compatible (gradient) fields, and both problems can be analyzed separately. The reaction term adds an interesting twist to the problem, coupling all five variational equations.	
\end{remark}

\begin{theorem}[Well-posedness of exact problem]
\label{theorem1}
	Problem P2 is well posed.
	
	\begin{proof}
		Following the theory of mixed problems \cite{brezzi-fortin,ern-guermond} one needs to show that $a(\cdot,\cdot)$ and $b(\cdot,\cdot)$ are continuous, that $a(\cdot,\cdot)$ is coercive in $\delta Z$, with
		\begin{eqnarray}
			\delta Z&=&\{x\in X~|~b(x,y)=0\, \quad \forall y\in Y\}\nonumber \\
			&=&\{(\uP,\eP,\sP)\in X \mid \eP-\nabla \uP=0, ~\nabla \cdot \sP + \zeta \uP =0,~\left. ~\sP \cdot \mathbf{\check{n}}\right |_{\partial\Omega_N} = 0  \}~,
		\end{eqnarray}
		and that the inf-sup condition on $b(\cdot,\cdot)$,
		\begin{equation}
			\inf_{y\in Y}\sup_{x\in X} \frac{b(x,y)}{\|x\|_X\|y\|_Y} \geq k_b > 0~.
		\label{eqinfsup}
		\end{equation}
		holds true.
		
		The norms adopted in $X$ and $Y$ are: 
		\begin{eqnarray}
			\|x\|_X &=& \|(\uP,\eP,\sP)\|_X = \sqrt{\|\nabla \uP\|^2+\|\eP\|^2+\kappa \|\sP\|^2}~,\\
			\|y\|_Y &=& \|(\lP,\mP)\|_Y=\sqrt{\frac{1}{\kappa}\left \| \nabla \lP \right \|^2+\left\|\mP \right \|^2}~.
		\end{eqnarray}
		
		They were chosen so that they have the same units. It is easy to see that both $X$ and $Y$ are Hilbert spaces, and that $a(\cdot,\cdot)$ and $b(\cdot,\cdot)$ are continuous. Notice that
		\begin{equation*}
				a(x,x)=\|\eP\|^2 + \kappa \|\sP\|^2 ~.
		\end{equation*}
		So, for $x \in \delta Z$ it holds that:
		\begin{equation}
			a(x,x)= \frac12 \|\eP\|^2 +\frac12 \|\nabla \uP\|^2+ \kappa \|\sP\|^2 \geq \frac12 \|x\|_X^2~.
		\label{eqcoercivity}
		\end{equation}
		and, thus, coercivity holds with constant equal to $\frac12$.
		
		It remains to show that there exists $k_b>0$ such that:
		\begin{equation*}
			\sup_{(\uP,\eP,\sP)\in X} \frac{\left (\mP,\eP-\nabla \uP \right ) - \left (\nabla \lP, \sP \right ) + \left ( \lP, \zeta \uP \right )}{\|(\uP,\eP,\sP)\|_X} \geq k_b \, \|(\lP,\mP)\|_{Y}~, \qquad \forall (\lP,\mP)\in Y~.
		\end{equation*}
		To this end, we can take $\uP=0$, $\eP=\mP$ and $\sP=-\frac{1}{\kappa}\nabla \lP$ to get
		\begin{equation*}
			\sup_{(\uP,\eP,\sP)\in X} \frac{\left (\mP,\eP-\nabla \uP \right ) - \left (\nabla \lP, \sP \right ) + \left ( \lP, \zeta\,\uP \right )}{\|(\uP,\eP,\sP)\|_X} 
			\geq \frac{\|\mP\|^2 + \frac{1}{\kappa}\|\nabla \lP\|^2}{\sqrt{\|\mP\|^2 + \frac{1}{\kappa}\|\nabla \lP \|^2}} = \|(\lP,\mP)\|_Y~,
		\end{equation*}
		giving an inf-sup constant $k_b$ equal to one. Thus, the coercivity of $a(\cdot,\cdot)$ on $\delta Z$ and the inf-sup stability of $b(\cdot,\cdot)$ have been proven. Then we have the well-posedness of the exact problem, with a unique solution $(\uP,\eP,\sP)\in H^1(\Omega) \times \left[L^2(\Omega)\right]^{d_{\Omega}} \times L^2(\Omega)$ that depends continuously on $q\in L^2(\Omega)$, $\etP \in \left[L^2(\Omega)\right]^{d_{\Omega}}$ and $\stP \in \left[L^2(\Omega)\right]^{d_{\Omega}}$.
	\end{proof}
\end{theorem}

\subsection{Approximation}

Let $U_h$, $E_h$, $S_h$, $L_h$ and $M_h$ be the finite element spaces for the five unknowns. We assume that $U_h\subset U$ and $L_h\subset L$ (in particular, the Dirichlet boundary conditions are strongly enforced). 
Also, we assume that $E_h$, $S_h$ and $M_h$ are subspaces of $[L^2(\Omega)]^d_{\Omega}$ and define $X_h:=U_h\times E_h \times S_h$ and $Y_h:=L_h\times M_h$. The numerical approximation we propose corresponds to solving the following discrete problem:

\begin{tcolorbox}[colback=black!10!white, boxrule=0pt, halign=justify]
	\begin{problem*}[P2$_{h}$]
		Find $(x_h,y_h)=\left ( \uP_h, \eP_h, \sP_h, \lP_h, \mP_h\right ) \in X_h\times Y_h$ such that:
		
		\begin{equation}
			(x_h,y_h)= \arg \left \{ \inf_{\bar{x}\in X_h} \sup_{\bar{y}\in Y_h} 
			\tilde{\mathcal{L}}(\bar{x},\bar{y})
			\right \}
		\end{equation}
		where
		\begin{eqnarray}
			\tilde{\mathcal{L}}(\bar{x},\bar{y})&=&J(\bar{\uP},\bar{\eP},\bar{\sP}) \nonumber \\
			&-&\left ( \bar{\sP},\nabla \bar{\lP} \right ) + \left (\zeta \bar{\uP} - q , \bar{\lP} \right )  
			+ \left ( \bar{\eP}-\nabla \bar{\uP}, \bar{\mP} \right ) \nonumber \\
			&+&\frac{\alpha}{2} \left \| \bar{\eP} - \nabla \bar{\uP} \right \|^2 \nonumber \\
			&-&\frac{\gamma}{2} \|\bar{\eP}+\bar{\mP}-\etP\|^2 \nonumber \\
			&-&\frac{\eta}{2\kappa} \|\kappa \bar{\sP} - \nabla \bar{\lP} - \kappa \stP\|^2 \nonumber \\
			&+& \frac{\theta \kappa \ell_{\sP}^2 }{2}\, \|\nabla \cdot \bar{\sP} + \zeta\, \bar{\uP} - q \|_h^2 \nonumber \\
			&-& \frac{\beta \ell_{\mP}^2 }{2}\, \|\nabla \cdot \bar{\mP} + \zeta\, \bar{\lP} \|_h^2~. \label{eq:discLag}
		\end{eqnarray}
	\end{problem*}
\end{tcolorbox}
\noindent
in which we kept the saddle-point structure by adding the squared norms of the optimality conditions for the dual variables and subtracting the square norms of the optimality conditions for the primal ones.


The mesh-dependent norm $\|\cdot\|_h$ is defined as 
$\|f\|_h^2 = \sum_K \|f\|^2_{L^2(K)} $
with the sum performed over all elements $K$ of the mesh. The associated scalar product is denoted as $\left ( \cdot, \cdot \right )_h$. This mesh-dependent norm is only necessary when the finite element spaces for $\mathbf{s}$ or $\boldsymbol{\mu}$ are discontinuous, since their divergence is then not in $L^2(\Omega)$. For continuous approximations of $\mathbf{s}$ and $\boldsymbol{\mu}$ the norm $\|\cdot\|_h$ coincides with the $L^2(\Omega)$-norm.

Above, $\alpha$, $\gamma$, $\eta$, $\theta$ and $\beta$ are non-negative non-dimensional coefficients to be further discussed shortly, and $\ell_{\sP}$ and $\ell_{\mP}$ are length scales. The function $\tilde{\mathcal{L}}$ is an augmented Lagrangian, of which the saddle point is the solution of P2$_h$. 

The optimality conditions of P2$_h$ are:

\begin{eqnarray}
	\alpha \left ( \nabla \uP_h, \nabla \delta \uP \right ) + \theta \kappa \ell_{\sP}^2  \left ( \zeta^2 \uP_h, \delta \uP \right ) - \alpha \left ( \eP_h,\nabla \delta \uP \right ) + & & \nonumber \\ 
	 + \theta \kappa \ell_f^2  \left ( \nabla \cdot \sP_h, \zeta \delta \uP \right )_h + \left ( \lP_h, \zeta\delta \uP \right ) - \left ( \mP_h, \nabla \delta \uP \right ) & = & \theta \kappa \ell_{\sP}^2  \left ( q, \zeta\delta \uP \right )\nonumber \\ & & \quad \forall \,\delta \uP \in U_h  \label{eq:var1h} \\ 
	-\alpha \left ( \nabla \uP_h, \delta \eP \right ) +(1+\alpha-\gamma)\,\left ( \eP_h,\delta \eP \right )+ (1-\gamma)\, \left ( \mP_h, \delta \eP \right ) &=& \left (1-\gamma)\,( \etP, \delta \eP \right )  \nonumber \\ & &\quad \forall \,\delta \eP \in E_h\label{eq:var2h}\\
	(1-\eta)\kappa \left ( \sP_h, \delta \sP \right ) - (1-\eta)\left ( \nabla \lP_h, \delta \sP \right ) + & & \nonumber \\
	+\theta \kappa \ell_{\sP}^2  \left ( \nabla \cdot \sP_h, \nabla \cdot \delta \sP \right )_h+ \theta \kappa \ell_{\sP}^2 \,\left (\zeta\,\uP_h, \nabla \cdot \delta \sP \right )_h & = & (1-\eta)\kappa \left ( \stP,\delta \sP \right ) + \nonumber \\ & &+\theta \kappa \ell_{\sP}^2  \left ( q , \nabla \cdot \delta \sP \right )_h \nonumber \\ & & \quad \forall\, \delta \sP \in S_h \label{eq:var3h}\\
	-(1-\eta)\left ( \sP_h,\nabla \delta \lP\right ) + \left ( \zeta\,\uP_h,\delta \lP \right ) - \beta \ell_{\mP}^2 \left ( \zeta^2 \lP_h, \delta \lP\right )+ & & \nonumber \\
	- \beta \ell_{\mP}^2 \left ( \nabla \cdot \mP_h, \zeta\,\delta \lP \right )_h -\frac{\eta}{\kappa} \left ( \nabla \lP_h , \nabla \delta \lP \right )&=&\left ( q,\delta \lP \right ) + \eta \left ( \stP,\nabla \delta \lP \right ) \nonumber \\ & & \quad \forall\,\delta \lP \in L_h\label{eq:var4h} \\
	-\left (\nabla \uP_h,\delta \mP \right )+(1-\gamma)\left ( \eP_h,\delta \mP \right )
	- \beta \ell_{\mP}^2 \left ( \zeta\,\lP_h, \nabla \cdot \delta \mP \right )_h + & & \nonumber \\- \beta \ell_{\mP}^2 \left ( \nabla \cdot \mP_h, \nabla \cdot \delta \mP \right )_h - \gamma \left ( \mP_h,\delta \mP \right ) & = & -\gamma \left (\etP,\delta \mP \right ) \nonumber \\ & & \quad \forall\,\delta \mP \in M_h \label{eq:var5h}
\end{eqnarray}

The sum of the left-hand sides (with the last two multiplied by $-1$) constitutes the discrete bilinear form that corresponds to P2$_h$, i.e.,
\begin{equation}
\tilde{B}((x_h,y_h),(\delta x,\delta y))=\tilde{a}(x_h,\delta x)+\tilde{b}(\delta x,y_h)-\tilde{b}(x_h,\delta y)+\tilde{c}(y_h,\delta y)~,
\end{equation}
where
\begin{equation}
	\begin{array}{rcl}
		\tilde{a}(x_h,\delta x)~&=&~(1-\gamma)\left (\eP_h,\delta \eP\right )+(1-\eta)\kappa \left ( \sP_h,\delta \sP \right ) +\\
		~&+&~\alpha \left( \eP_h-\nabla \uP_h, \delta \eP - \nabla \delta \uP \right)+\\
		~&+&~ \theta \kappa \ell_{\sP}^2  \left ( \nabla \cdot \sP_h + \zeta \uP_h , \nabla \cdot \delta \sP + \zeta \, \delta \uP \right )_h~,
	\end{array}
\end{equation}

\begin{equation}
	\begin{aligned}
		\tilde{b}(\delta x,y_h)~=~b(\delta x,y_h)-\gamma \left (\mP_h,\delta \eP\right )+\eta \left (\nabla \lP_h,\delta \sP \right ) ~,
	\end{aligned}
\end{equation}

\begin{equation}
	\begin{array}{rcl}
		\tilde{c}(y_h,\delta y) ~&=&~  \beta \ell_{\mP}^2 \left ( \nabla \cdot \mP_h + \zeta \lP_h , \nabla \cdot \delta \mP + \zeta\,\delta \lP \right )_h +\\
		~&+&~\dfrac{\eta}{\kappa}\,\left ( \nabla \lP_h,\nabla \delta \lP \right ) + \gamma \left ( \mP_h,\delta \mP \right )~,
	\end{array}
\end{equation}
and $b(\cdot,\cdot)$ as defined in (\ref{eq:bform}).

Problem P2$_h$ can thus be written as:

\begin{tcolorbox}[colback=black!10!white, boxrule=0pt, halign=justify]
	Find $(x_h,y_h)\in X_h \times Y_h$ such that:
	\begin{equation}
		\tilde{B}\left ( \left ( x_h,y_h \right ),\left ( \delta x,\delta y \right ) \right )= \tilde{F}(\delta x,\delta y), \quad \forall (\delta x, \delta y)\in X_h\times Y_h
		\label{eq:BzF}
	\end{equation}
\end{tcolorbox}

The linear form $\tilde{F}:X_h\times Y_h \to \mathbb{R}$ stands for the sum of the right-hand sides of (\ref{eq:var1h})-(\ref{eq:var5h}), with the last two multiplied by $-1$.

We now define the formulations considered in the analysis:

\begin{description}
\item[(a)] {\bf Natural formulation}. It corresponds to a choice of discrete spaces that satisfies
$$
\nabla U_h \subseteq M_h, \qquad \nabla L_h \subseteq S_h, \qquad E_h = M_h 
$$
When such is the case, as we will see, no stabilization is required. The coefficients are thus taken equal to zero,
$$
\alpha = \gamma = \eta = \theta = \beta = 0~.
$$
The typical spaces for this formulation are continuous piecewise $P_{k+1}$ for $U_h$ and $L_h$, and discontinuous, piecewise $P_{k}$ for $M_h$, $S_h$ and $E_h$, $k \ge 0$.

\item[(b)] {\bf Equal-order formulation}. It corresponds to choosing the same $H^1$-conforming finite element space for all unknown fields, typically a continuous piecewise $P_k$ approximation. Equal-order formulations require stabilization. We consider three choices of the stabilization parameters:

\begin{description}

\item[(b-1)] {\bf Unstabilized method:}
All coefficients set to zero. This method is not optimally convergent.

\item[(b-2)] {\bf Minimally-stabilized method:} Only coefficients $\alpha$ and $\eta$ are positive. In particular, we take
$$
\alpha = \frac{1}{8}, \qquad \eta = \frac12, \qquad \gamma = \theta = \beta = 0~.
$$

\item[(b-3)] {\bf Fully-stabilized method:} All coefficients are positive, and the scale lengths are equal to the mesh size; i.e.,
	\begin{equation}
	\alpha = \frac{1}{8}, 
	\quad \gamma=\alpha, \quad \eta = \frac12,\quad \theta=\frac12,\quad \beta=\frac12,
	\quad \ell_{\sP} = \ell_{\mP} = h~. \label{eq:paramsopt}
	\end{equation}

\end{description}

\end{description} 

\noindent With these definitions, it is possible to prove

\bigskip

\begin{theorem}[Well-posedness of discrete problem and optimal approximation] For methods (a), (b-2) and (b-3) above, the discrete problem P2$_h$ is well-posed, and $(x_h,y_h)$ satisfies that $\exists c>0$ independent of $h$ such that
	\begin{equation}
		\||(x,y)-(x_h,y_h)\|| \leq c\,\min_{(\bar{x},\bar{y})\in X_h\times Y_h}
		\||(x,y)-(\bar{x},\bar{y})\||
		\label{eq:bestapprox}
	\end{equation}
where the norm $\||\cdot\||$ depends on the method: For methods (a) and (b-2) it is
\begin{equation}
\||(\bar{x},\bar{y})\||^2 ~=~ \|(\bar{x},\bar{y})\|^2_{X\times Y} ~=~ \|\nabla \bar{\uP} \|^2 +\|\bar{\eP}\|^2 + \kappa \|\bar{\sP}\|^2 +  \frac{1}{\kappa} \|\nabla \bar{\lP}\|^2 + \|\bar{\mP}\|^2 ~,
\end{equation}
and for method (b-3) it is 
\begin{equation}
\||(\bar{x},\bar{y})\||^2 ~=~ \|(\bar{x},\bar{y})\|^2_{X\times Y}  + \kappa h^2 \|\nabla \cdot \bar{\sP}\|^2  + h^2\,\|\nabla \cdot \bar{\mP}\|^2~.
	\end{equation}
\end{theorem}


\begin{proof}
	It is easy to verify that the exact solution $(x,y)$ satisfies the discrete problem (\ref{eq:BzF}) when substituted in the place of $(x_h,y_h)$. The formulation is thus consistent. The continuity of all the involved bilinear and linear forms in the norm $\||\cdot\||$ is straightforward. The rest of the proof is different for each method.
	
	\bigskip
	
 \noindent{\bf (a) Natural formulation}. It is the classical mixed approximation, in which the discrete bilinear forms are $\tilde{a}=a$ and $\tilde{b}=b$ and $\tilde{c}=0$. Let
  $$
 \delta Z_h=\{ x_h\in X_h ~|~b(x_h,y_h)=0\quad \forall y_h\in Y_h \}
 $$
be the discrete kernel of $b(\cdot,\cdot)$. Since $\nabla U_h\subset E_h=M_h$ we have $\eP_h = \nabla \uP_h$ and thus $\delta Z_h \subset \delta Z$. We thus have coercivity in the discrete kernel ($a(x_h,x_h) \geq \frac12 \|x_h\|_X^2$ for all $x_h \in \delta Z_h$) immediately from (\ref{eqcoercivity}). 

According to Brezzi's theorem \cite{Demkowicz2006} it only remains to show that the inf-sup condition (\ref{eqinfsup}) holds for the discrete spaces $X_h$ and $Y_h$ with a mesh-independent constant. Consider an arbitrary $y_h = \left ( \lP_h, \mP_h \right )\in Y_h$ and take $\uP_h=0$, $\eP_h=\mP_h$ and $\sP_h = -\frac{1}{\kappa} \nabla \lP_h$. These two latter choices are possible, since $E_h=M_h$ and $\nabla L_h \subset S_h$, which yields
 \begin{equation}
 \sup_{x_h=(\uP_h,\eP_h,\sP_h)\in X_h} \frac{b\left (x_h,y_h\right )}{\|x_h\|_X} \geq \|y_h\|_Y~,
 \label{eqsupb}
 \end{equation}
 with inf-sup constant equal to 1.
 
 \bigskip
 
 \noindent{\bf (b-2) Equal-order minimally-stabilized method.} For this case we use Babuska's theorem \cite{Demkowicz2006}, which requires proving that $\exists K_B>0$, independent of $h$, such that
 \begin{equation}
 \sup_{(\bar{x},\bar{y})\in X_h\times Y_h} \frac{\tilde{B}\left ( (x_h,y_h), (\bar{x},\bar{y}) \right )}
 {\||(\bar{x},\bar{y})\||} \geq K_B\,\||(x_h,y_h)\||~,
 \end{equation}
 for all $(x_h,y_h)$ in $X_h\times Y_h$. Let $k_0 > 0$ and
 $(\bar{x},\bar{y})=\left (\uP_h, \eP_h - k_0 \mP_h, \sP_h, \lP_h, \mP_h \right )$. Then,
 \begin{eqnarray}
 \tilde{B}\left ((x_h,y_h), (\bar{x},\bar{y}) \right ) & = &
 \alpha \|\nabla \uP_h\|^2 + (1+\alpha) \|\eP_h\|^2 + (1-\eta) \kappa \|\sP_h\|^2 + \nonumber \\
 & & + k_0 \|\mP_h\|^2 + k_0 \alpha \left ( \nabla \uP_h, \mP_h \right ) 
 - k_0 (1+\alpha) \left ( \eP_h,\mP_h \right ) + \nonumber \\
 & & - 2\alpha \left ( \eP_h , \nabla \uP_h \right ) +\frac{\eta}{\kappa} \|\nabla \lP_h\|^2
 + \eta \left ( \eP_h , \mP_h \right )~.
 \end{eqnarray}
 Now we use Young's inequality ($2ab\geq -k a^2-(1/k)b^2$, $\forall a, b \in \mathbb{R}, \forall k>0$) as follows
 \begin{eqnarray*}
 -2\alpha \left ( \eP_h, \nabla \uP_h \right ) & \geq &
 -\frac{1}{k_1} \| \eP_h \|^2 - k_1 \alpha^2 \|\nabla \uP_h \|^2 \\
 \left ( \nabla \uP_h, \mP_h \right ) & \geq & -\frac{1}{2k_2} 
 - \frac{k_2}{2} \|\mP_h\|^2 \\
 -\left ( \eP_h, \mP_h \right ) & \geq & -\frac{1}{2 k_3} \|\eP_h \|^2
 - \frac{k_3}{2} \|\mP_h\|^2 \\
 \eta \left ( \eP_h, \mP_h \right ) & \geq & -\frac{\eta}{2 k_4} \|\mP_h\|^2
 - \frac{\eta k_4}{2} \|\eP_h\|^2
 \end{eqnarray*}
 to arrive at
 \begin{eqnarray*}
 \tilde{B}\left ((x_h,y_h), (\bar{x},\bar{y}) \right ) & \geq & 
 A_1 \|\nabla \uP_h\|^2 + A_2 \|\eP_h\|^2 + (1-\eta) \kappa \|\sP_h\|^2 +
 \frac{\eta}{\kappa} \|\lP_h\|^2 + A_3 \|\mP_h\|^2~,
 \end{eqnarray*}
 with
 \begin{eqnarray*}
 A_1 & = & \alpha - k_1 \alpha^2 - \frac{\alpha k_0}{2 k_2}~, \\
 A_2 & = & 1+\alpha - \frac{1}{k_1} - \frac{(1+\alpha)k_0}{2 k_3} - \frac{\eta k_4}{2}~, \\
 A_3 & = & k_0 - \frac{\alpha k_0 k_2}{2} - \frac{(1+\alpha)k_0k_3}{2} - \frac{\eta}{2 k_4}~.
 \end{eqnarray*}
 The above expressions are valid for any $k_i>0$, $i=0,\ldots,4$. Numerical optimization was used to determine the following set of parameters
 $$
 \alpha=\frac18,\quad \eta=\frac12,\quad k_0=0.352, \quad k_1=2.606, \quad
 k_2=0.606,\quad k_3=0.707, \quad k_4=1.654
 $$
 which produce strictly positive values $A_1 \simeq A_2 \simeq A_3 \simeq 0.047$.  At the same time, it is not difficult to prove that, when $k_0=0.352$, it holds that $\||(\bar{x},\bar{y})\||\leq \sqrt{2}\,\||(x_h,y_h)\||$, which leads to
 $$
 K_B \geq \frac{0.047}{\sqrt{2}} \simeq 0.033~,
 $$
 independent of $h$ as required.
 
 \bigskip
 
 \noindent{\bf (b-3) Equal-order fully-stabilized method.} Similarly to (b-2), this case is proved using Babuska's theorem \cite{Demkowicz2006}. There is a difference, however, since when all stabilizing terms are active the bilinear form $\tilde{B}(\cdot,\cdot)$ is {\em strongly} coercive in the norm $\||\cdot\||$, i.e., $\exists K_B > 0$ (independent of $h$) such that
	\begin{equation}
		\tilde{B}((x_h,y_h),(x_h,y_h))\geq K_B \left ( 
		\|x_h\|_X^2 + \|y_h\|_Y^2 + \kappa\,h^2 \|\nabla \cdot \sP_h \|^2
		+ h^2\,\|\nabla \cdot \mP_h \|^2\right ),
	\end{equation}
	for all $ (x_h,y_h)\in X_h\times Y_h$~.
	By direct computation, 
	\begin{eqnarray}
		\tilde{B}((x_h,y_h),(x_h,y_h)) ~&=&~ (1-\gamma) \|\eP_h\|^2 + (1-\eta)\kappa \|\sP_h\|^2 \nonumber \\
		~&+&~ \alpha \|\eP_h-\nabla \uP_h\|^2 + \theta \kappa \ell_{\sP}^2 \|\nabla \cdot \sP_h + \zeta \uP_h\|^2 \nonumber \\
		~&+&~ \beta \ell_{\mP}^2 \| \nabla \cdot \mP_h + \zeta \lP_h \|^2
		+ \frac{\eta}{\kappa} \|\nabla \lP_h\|^2 + \gamma \|\mP_h\|^2~.
	\end{eqnarray}
	
	Now, as before, we apply Young's inequality to each of the cross-product terms in the previous expression.
	
	\begin{eqnarray}
		-2\alpha \left (\eP_h,\nabla \uP_h \right ) & \geq & -\frac14 \|\eP_h\|^2 - 4\alpha^2 \|\nabla \uP_h \|^2 \\
		2\theta\kappa \ell_{\sP}^2\left ( \nabla \cdot \sP_h, \zeta \uP_h \right )_h& \geq & - \frac12 \kappa \theta \ell_{\sP}^2 \|\nabla\cdot \sP_h\|^2 - {2\theta\kappa\ell_{\sP}^2}\|\zeta \uP_h\|^2\\
		2 \beta \ell_{\mP}^2\left ( \nabla \cdot \mP_h, \zeta \lP_h \right )_h& \geq & - \frac12 {\beta} \ell_{\mP}^2 \|\nabla\cdot \mP_h\|^2 - {2\beta\ell_{\mP}^2}\|\zeta \lP_h\|^2
	\end{eqnarray}
	
	Then, adding up all terms that contain $\sP_h$ we arrive at
	\begin{equation}
		(1-\eta) \kappa \|\sP_h\|^2 + \frac12 \kappa \ell_{\sP}^2 
		\theta  \|\nabla\cdot \sP_h\|_h^2~.
		\label{eq:term1}
	\end{equation}
	Adding up all terms that contain $\zeta \uP_h$ gives
	$$
	-\theta \kappa \ell_{\sP}^2 \|\zeta \uP_h\|^2
	$$
	which using Poincar\'e-Friedrichs inequality ($\|\uP_h\| \leq c_{\mbox{\scriptsize{PF}}}\,\ell_{\Omega}\,\|\nabla \uP_h\|$, 
	with $\ell_{\Omega}$ the diameter of $\Omega$) is shown to be greater than or equal to 
	\begin{equation}
		- {\theta \kappa \ell_{\sP}^2 \|\zeta\|_{\infty}^2 c_{\mbox{\scriptsize{PF}}}^2 \ell_{\Omega}^2}\,\|\nabla \uP_h \|^2~.
		\label{eq:term2}
	\end{equation}
	By $\|\zeta\|_{\infty}$ we denote the $L^{\infty}(\Omega)$-norm of $\zeta$. Adding up all terms that contain $\mP_h$ we arrive at
	\begin{equation}
		\gamma \|\mP_h\|^2 +\frac12 \beta\ell_{\mP}^2  \|\nabla \cdot \mP_h \|_h^2~.
		\label{eq:term4}
	\end{equation}
	The sum of the terms that contain $\zeta \lP_h$ is, by Poincaré-Friedrichs inequality, greater than or equal to
	\begin{equation}
		- \beta \ell_{\mP}^2 \|\zeta\|_{\infty}^2 c_{\mbox{\scriptsize{PF}}}^2 \ell_{\Omega}^2\,\|\nabla \lP_h \|^2~.
		\label{eq:term5}
	\end{equation}
	Adding up all terms that contain $\nabla \uP_h$ gives
	\begin{equation}
		\left ( \alpha - 4 \alpha^2 \right )
		\|\nabla \uP_h\|^2~.
		\label{eq:term6}
	\end{equation}
	Adding up all terms that contain $\eP_h$ gives
	\begin{equation}
		\left ( \frac34-\gamma + \alpha \right ) \|\eP_h\|^2~.
		\label{eq:term7}
	\end{equation}
	Adding everything up we get
	\begin{eqnarray}
		\tilde{B}((x_h,y_h),(x_h,y_h))~&\geq&~ \left (\frac34 -\gamma+\alpha \right ) \|\eP_h\|^2 + (1-\eta)\kappa \|\sP_h\|^2 + \nonumber \\
		~&+&~ (\alpha - 4 \alpha^2 - \theta \kappa \ell_{\sP}^2 \|\zeta\|_{\infty}^2 c_{\mbox{\scriptsize{PF}}}^2 \ell_{\Omega}^2)\, \|\nabla \uP_h\|^2 + \nonumber \\
		~&+&~\frac12 \theta \kappa \ell_{\sP}^2 \|\nabla \cdot \sP_h\|^2 + \nonumber \\
		~&+&~\frac12 \beta \ell_{\mP}^2 \| \nabla \cdot \mP_h\|^2 + \nonumber \\
		~&+&~\left ( \frac{\eta}{\kappa}-\beta \ell_{\mP}^2 \|\zeta\|_{\infty}^2 c_{\mbox{\scriptsize{PF}}}^2 \ell_{\Omega}^2 \right ) \|\nabla \lP_h\|^2 + \gamma \|\mP_h\|^2~. \label{eq:Bhparams}
	\end{eqnarray}
	Taking
	\begin{equation}
	\alpha = \frac{1}{8}, \quad \theta=\frac12,\quad \beta=\frac12,
	\quad \gamma=\alpha, \quad \eta = \frac12,
	\quad \ell_{\sP} = \ell_{\mP} = h, \label{eq:paramsopt}
	\end{equation}
	we get
	\begin{eqnarray*}
		\tilde{B}\left ( (x_h,y_h),(x_h,y_h) \right ) &\geq  &
		\left ( \frac{1}{16} - O(h^2) \right )\, \|\nabla \uP_h\|^2+\\
		& & +\frac{3}{4}\,\|\eP_h\|^2+\frac12 \kappa \|\sP_h\|^2 + \frac14 \kappa h^2 \|\nabla \cdot \sP_h\|_h^2 + \\
		& &
		+\left (\frac{1}{2\kappa}  -O(h^2) \right ) \|\nabla \lP_h\|^2 + \\
		& & + \frac18 \|\mP_h\|^2  +\frac14 h^2 \|\nabla \cdot \mP_h\|_h^2 
	\end{eqnarray*}
	and finally
	\begin{equation}
		\tilde{B}\left ( (x_h,y_h),({x}_h,{y}_h) \right ) \geq 
		K_B \, \| |(x_h,y_h)  \| |_h^2 
	\end{equation}
	where
	\begin{equation}
		K_B \geq \frac{1}{16}~.
	\end{equation}
	We have assumed that $h$ is sufficiently small so that the terms $O(h^2)$ can be neglected and thus $K_B$ becomes mesh independent. 
	
	\bigskip
	
	In all three cases we have proved that $\tilde{B}$ is consistent, continuous and coercive, which implies that the discrete problem is well posed and that the best approximation property (\ref{eq:bestapprox}) holds true.
\end{proof}

\section{Numerical experiments}

In this section, we numerically assess the convergence properties of the finite element approximation to the proposed problem. We consider a sequence of finite element conformal and shape-regular partitions of the computational domain $\mathcal{T}_h$ with characteristic mesh size $h$. 
Let us define the discrete spaces associated to $\mathcal{T}_h$
$$
\CG_{k}(\mathcal{T}_h) = \{v \in C^0(\Omega), v|_K \in \mathbb{P}_k(K)~\forall K \in \mathcal{T}_h\}
$$
and 
$$
\DG_{k}(\mathcal{T}_h) = \{v, v|_K \in \mathbb{P}_k(K)~\forall K \in \mathcal{T}_h\}
$$
where $\mathbb{P}_k(K)$ is the space of polynomials of degree $k$ on $K$.

Two different choices of finite element spaces are then considered to approximate
the optimization problem, namely, a \textbf{natural formulation}, in which we take 
\begin{equation}\label{eq:natural_func_spaces}
	\left(\uP_{h}, \eP_{h}, \sP_{h}, \lP_{h}, \mP_{h}\right) \in \CG_{k+1} \times \left[\DG_{k}\right]^{d_{\Omega}} \times \left[\DG_{k}\right]^{d_{\Omega}} \times \CG_{k+1} \times \left[\DG_{k}\right]^{d_{\Omega}} 
\end{equation} 
and an \textbf{equal-order formulation}, in which we take
\begin{equation}\label{eq:equal_order_func_spaces}
\left(\uP_{h}, \eP_{h}, \sP_{h}, \lP_{h}, \mP_{h}\right)  \in \CG_{k+1} \times \left[\CG_{k+1}\right]^{d_{\Omega}} \times \left[\CG_{k+1}\right]^{d_{\Omega}} \times \CG_{k+1} \times \left[\CG_{k+1}\right]^{d_{\Omega}}
\end{equation}
In the numerical tests we take $k \in \{0, 1, 2\}$ for both formulations, corresponding to the use of linear, quadratic or cubic polynomials for the scalar fields $\uP$ and $\lP$ in both formulations. Indeed, the Natural formulation
considers one degree lower for the vector fields $\sP$, $\eP$ and $\mP$, whereas for the equal-order formulation, the degree is the same for all fields. 
To assess the convergence properties of the method we consider both convex and non-convex geometries in two spatial dimensions. The convergence behavior is assessed by computing the error with respect to manufactured solutions.

\subsection{Test case 1: a convex domain}\label{sec:test_case1}

\begin{figure}[h!]
	\centering
	\caption{Figure (a) presents the functional relationship between the norm of the flux and the gradient $\|\sP\|~\mbox{vs.}~\|\eP\|$, and $\|\stP\|~\mbox{vs.}~\|\etP\|$, respectively. Figure (b) to (f) shows
		$\uP$, $\lP$, $\eP$, $\sP$, and $\mP$.}
	\def\svgwidth{1.0\columnwidth}
	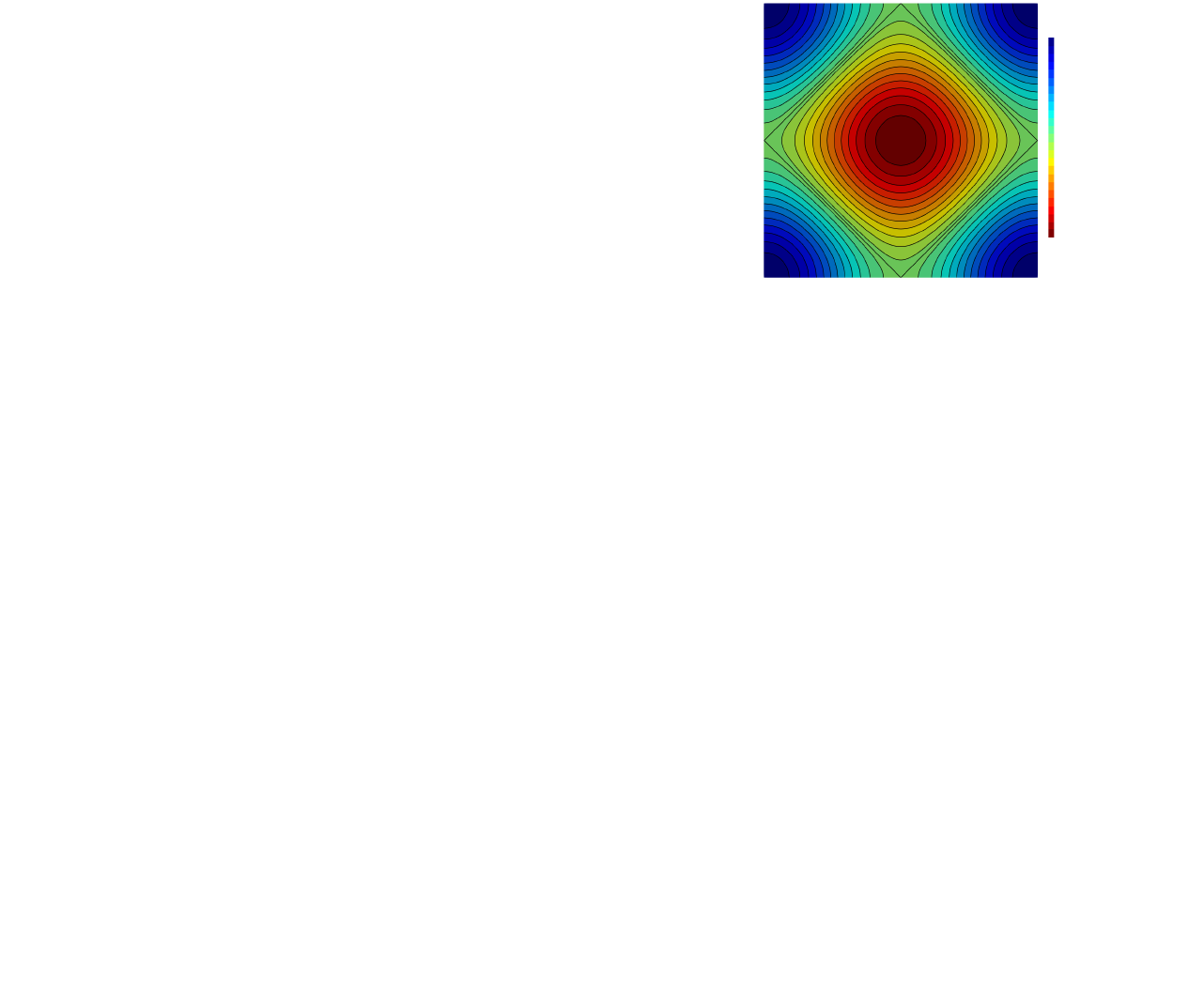
	\label{fig:all_fields_ref}
\end{figure}

Let us first describe the procedure for constructing the manufactured solution for this test problem. First, consider the smooth field $\uP$ defined by:
\begin{equation}
	\uP:\Omega \to \mathbb{R},\quad \uP(x,y) = \cos(\pi x)\cos(\pi y),
\end{equation}
whose gradient is given by:
\begin{equation*}
	\begin{aligned}
	\eP: \Omega \to \mathbb{R}^{2}, \quad \eP(x,y) = -\pi \, [\sin(\pi x)\cos(\pi y), \cos(\pi x)\sin(\pi y)]^{\top}. \\
	\end{aligned}
\end{equation*}
Also define 
$$
\lP: \Omega \rightarrow \mathbb{R},\qquad \lP(x,y) = -\frac{1}{40\,\pi}\,\left[\cos\left(2\,\pi\,x \right) + \cos\left(2\,\pi\,y \right)\right] \\ \\
$$
and
\begin{equation}
\sP:\Omega \rightarrow \mathbb{R}^2, \qquad \sP = -\dfrac{\eP}{\|\eP\|}\, \left(\|\eP\| - \dfrac{1}{40}\,\|\eP\|^{3}\right) \label{eq:constlaw}
\end{equation}
At this stage, solely for the purpose of assessing the convergence properties of the discrete formulations under general conditions, it is convenient to introduce an additional source term on the right-hand side of eq. (\ref{eq:dif1}), i.e.,
\begin{equation}
	\zeta \lP  + \nabla \cdot \mP  =  f \label{eq:balance_with_f}
\end{equation}
This amounts to adding the term $-({f , \uP})$ to the Lagrangian (\ref{eq:lagrangian}).
Notice that the addition of this linear term does not affect the convergence analysis from previous sections.
By examining the optimality conditions in their strong form (\ref{eq:dif1})-(\ref{eq:dif5}),
we proceed as follows to construct the manufactured solution for the remaining fields.
First take
$$
\etP:\Omega \rightarrow \mathbb{R}^2, \qquad \etP(x,y) = \eP(x,y) + \frac{1}{20} \left[\sin\left(4\,\pi\,x \right), \sin\left(4\,\pi\,y\right) \right]^{\top}
$$
as being a perturbation of $\eP$, from which $\mP$ can be obtained using eq. (\ref{eq:dif2}).
Now, $\stP$ is computed using eq. (\ref{eq:dif3}) and finally, the source 
terms $f$ and $q$ are obtained after substituting the proposed fields 
into eqs. (\ref{eq:dif1}) and (\ref{eq:dif4}). 
In the numerical experiments we take $\kappa = 1$ and $\zeta = 1$.
Notice that, in this way, the fields $\eP$ and $\sP$ satisfy a type of 
constitutive relationship, which according to (\ref{eq:constlaw}) is essentially
a linear law plus a cubic perturbation term, whereas the values taken by the fields 
$\etP$ and $\stP$ at different points throughout the domain do not necessarily do so. 
This mimics practical scenarios in which these are obtained from experimental 
measurements or synthetic data. This is illustrated in part (a) of figure \ref{fig:all_fields_ref} which shows
a plot of the norm of $\sP$ as a function of the norm of $\eP$, along with their counterparts $\stP$ and $\etP$ evaluated over a computational mesh. 
Notice that the exact flux-gradient relationship is depicted as a strictly monotonic curve (thick black line), whereas the fields $\stP$ and $\etP$ (blue dots) exhibit a certain dispersion around the former. Contours of the different fields are also shown in the
figure considering the computational domain $\Omega$ to be the unit square $[0,1]^2$. The boundary of $\Omega$ is defined as $\partial \Omega = \Gamma_{\text{left}} \cup \Gamma_{\text{right}} \cup \Gamma_{\text{bottom}} \cup \Gamma_{\text{top}}$.
The problem is solved imposing the following boundary conditions
\begin{equation}
	\left\lbrace
	\begin{array}{rclll}
		\sP \cdot \mathbf{\check{n}} 			&=& \mathsf{g}_{\sP} 			&\mbox{on}&~\Gamma_{\text{bottom}} \cup \Gamma_{\text{top}} \\
		\mP \cdot \mathbf{\check{n}} &=& \mathsf{g}_{\mP} &\mbox{on}&~\Gamma_{\text{bottom}} \cup \Gamma_{\text{top}} \\
		\uP &=& \mathsf{g}_{\uP} &\mbox{on}&~\Gamma_{\text{left}} \cup \Gamma_{\text{right}} \\
		\lP &=& \mathsf{g}_{\lP} &\mbox{on}&~\Gamma_{\text{left}} \cup \Gamma_{\text{right}} \\
	\end{array}
	\right.
\end{equation}
where the values in the right hand side are computed from the manufactured fields.

\begin{figure}[h!]
	\centering
	\caption{Convergence of scalar fields $\uP$ and $\lP$ with $k=0$. Figures~(a)--(b) show results for $\uP$, and (c)--(d) for $\lP$ in both Natural and equal-order formulations, evaluated in the $L^{2}$ and $H^{1}$ norms.   Three stabilization configurations are compared: Unstabilized, Fully stabilized, and Minimally stabilized ($\alpha$, $\eta$).
	}
	\def\svgwidth{1.0\columnwidth}
\begingroup%
  \makeatletter%
  \providecommand\color[2][]{%
    \errmessage{(Inkscape) Color is used for the text in Inkscape, but the package 'color.sty' is not loaded}%
    \renewcommand\color[2][]{}%
  }%
  \providecommand\transparent[1]{%
    \errmessage{(Inkscape) Transparency is used (non-zero) for the text in Inkscape, but the package 'transparent.sty' is not loaded}%
    \renewcommand\transparent[1]{}%
  }%
  \providecommand\rotatebox[2]{#2}%
  \newcommand*\fsize{\dimexpr\f@size pt\relax}%
  \newcommand*\lineheight[1]{\fontsize{\fsize}{#1\fsize}\selectfont}%
  \ifx\svgwidth\undefined%
    \setlength{\unitlength}{529.13387269bp}%
    \ifx\svgscale\undefined%
      \relax%
    \else%
      \setlength{\unitlength}{\unitlength * \real{\svgscale}}%
    \fi%
  \else%
    \setlength{\unitlength}{\svgwidth}%
  \fi%
  \global\let\svgwidth\undefined%
  \global\let\svgscale\undefined%
  \makeatother%
  \begin{picture}(1,0.8234623)%
    \lineheight{1}%
    \setlength\tabcolsep{0pt}%
    \put(0.24211432,0.42113335){\color[rgb]{0,0,0}\makebox(0,0)[lt]{\lineheight{1.25}\smash{\begin{tabular}[t]{l}$\mbox{(a)}$\end{tabular}}}}%
    \put(0.74382664,0.42102165){\color[rgb]{0,0,0}\makebox(0,0)[lt]{\lineheight{1.25}\smash{\begin{tabular}[t]{l}$\mbox{(b)}$\end{tabular}}}}%
    \put(0.24273353,0.00320637){\color[rgb]{0,0,0}\makebox(0,0)[lt]{\lineheight{1.25}\smash{\begin{tabular}[t]{l}$\mbox{(c)}$\end{tabular}}}}%
    \put(0.74382664,0.00320637){\color[rgb]{0,0,0}\makebox(0,0)[lt]{\lineheight{1.25}\smash{\begin{tabular}[t]{l}$\mbox{(d)}$\end{tabular}}}}%
    \put(0,0){\includegraphics[width=\unitlength,page=1]{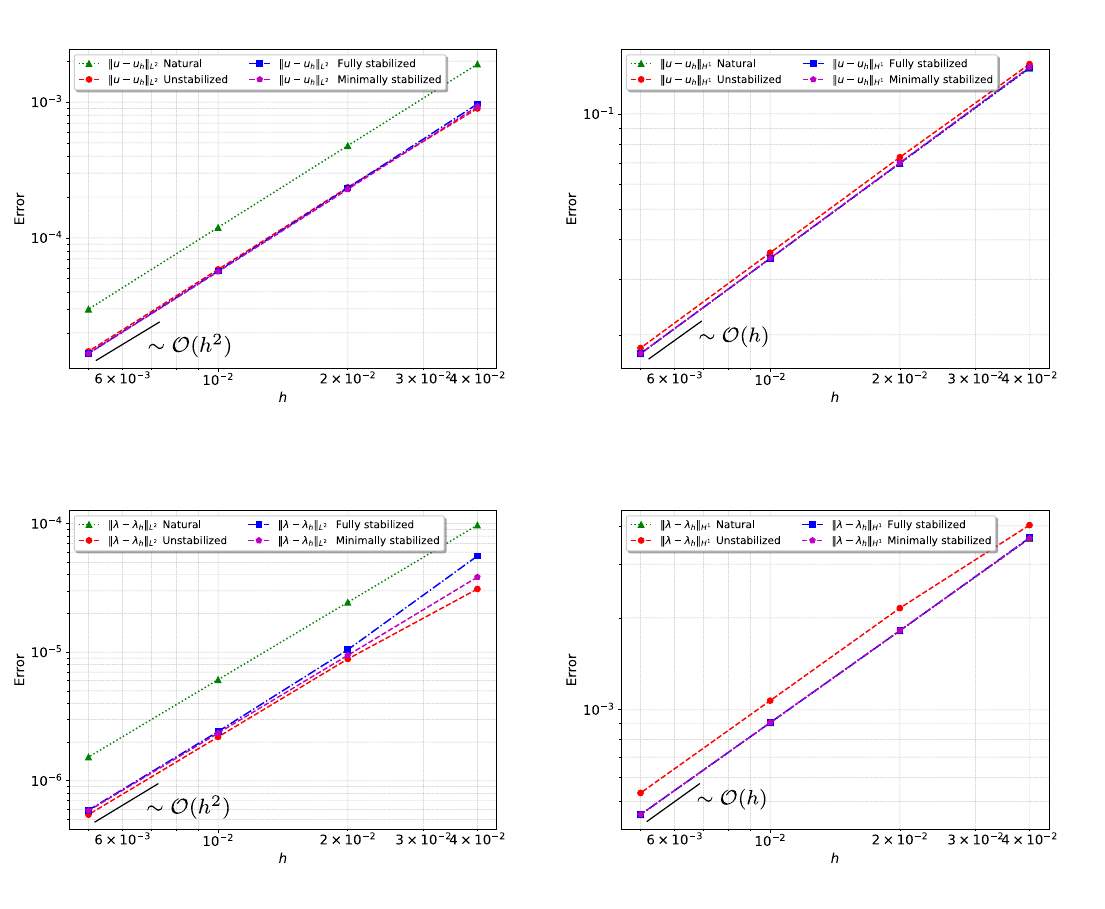}}%
  \end{picture}%
\endgroup%

	\label{fig:convergence_scalar_linear}
\end{figure}
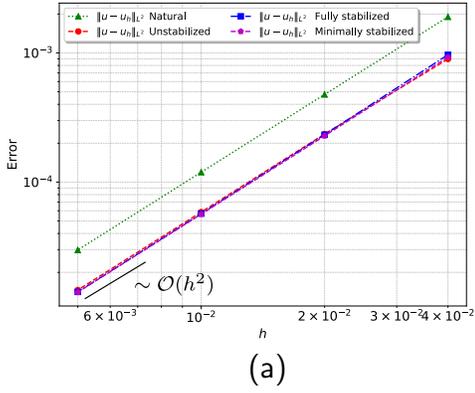

\begin{figure}[h!]
	\centering
	\caption{Convergence of vector fields $\sP$, $\eP$, and $\mP$ with $k=0$, evaluated in the $L^{2}$-norm for both Natural and equal-order formulations. Three stabilization configurations are compared: Unstabilized, Fully stabilized, and Minimally stabilized ($\alpha$, $\eta$).
	}
	\def\svgwidth{1.0\columnwidth}
\begingroup%
  \makeatletter%
  \providecommand\color[2][]{%
    \errmessage{(Inkscape) Color is used for the text in Inkscape, but the package 'color.sty' is not loaded}%
    \renewcommand\color[2][]{}%
  }%
  \providecommand\transparent[1]{%
    \errmessage{(Inkscape) Transparency is used (non-zero) for the text in Inkscape, but the package 'transparent.sty' is not loaded}%
    \renewcommand\transparent[1]{}%
  }%
  \providecommand\rotatebox[2]{#2}%
  \newcommand*\fsize{\dimexpr\f@size pt\relax}%
  \newcommand*\lineheight[1]{\fontsize{\fsize}{#1\fsize}\selectfont}%
  \ifx\svgwidth\undefined%
    \setlength{\unitlength}{529.13382943bp}%
    \ifx\svgscale\undefined%
      \relax%
    \else%
      \setlength{\unitlength}{\unitlength * \real{\svgscale}}%
    \fi%
  \else%
    \setlength{\unitlength}{\svgwidth}%
  \fi%
  \global\let\svgwidth\undefined%
  \global\let\svgscale\undefined%
  \makeatother%
  \begin{picture}(1,0.8220581)%
    \lineheight{1}%
    \setlength\tabcolsep{0pt}%
    \put(0,0){\includegraphics[width=\unitlength,page=1]{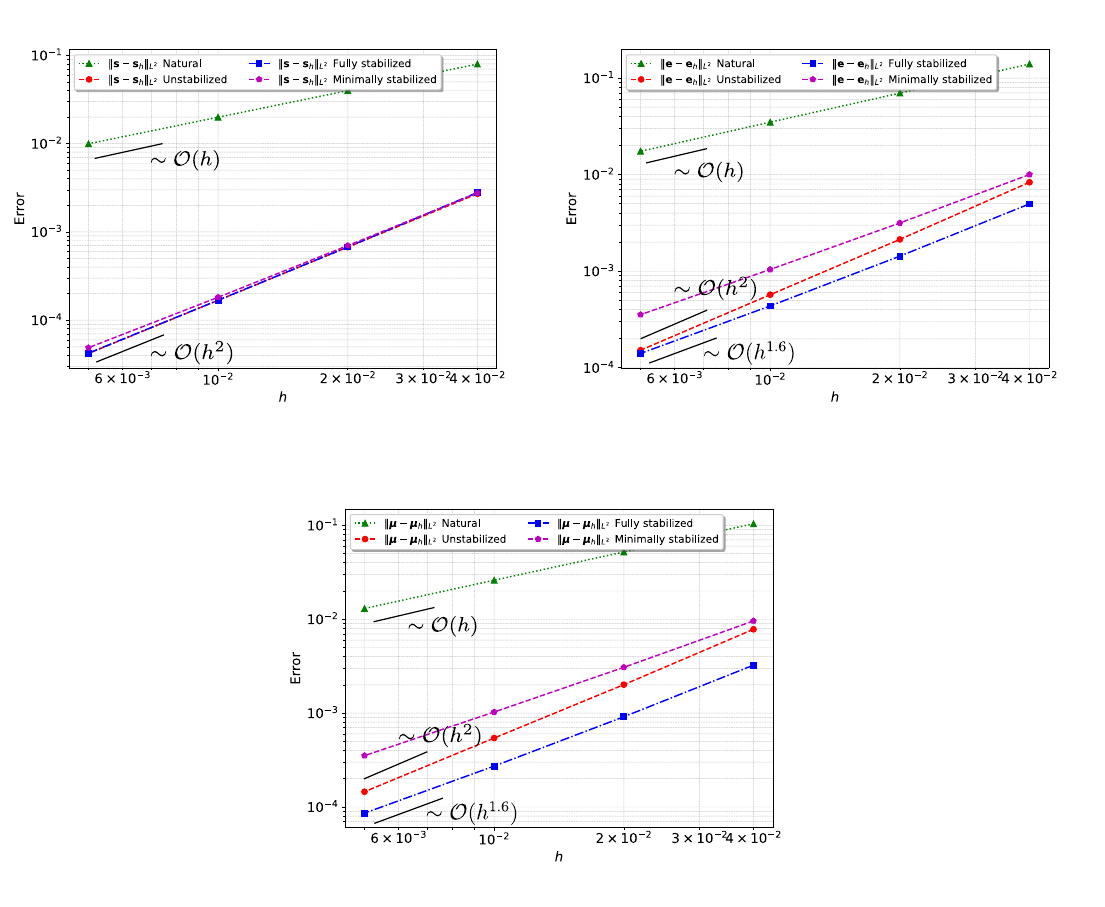}}%
    \put(0.24521589,0.42019719){\color[rgb]{0,0,0}\makebox(0,0)[lt]{\lineheight{1.25}\smash{\begin{tabular}[t]{l}$\mbox{(a)}$\end{tabular}}}}%
    \put(0.74976389,0.42019727){\color[rgb]{0,0,0}\makebox(0,0)[lt]{\lineheight{1.25}\smash{\begin{tabular}[t]{l}$\mbox{(b)}$\end{tabular}}}}%
    \put(0.4939187,0.00320651){\color[rgb]{0,0,0}\makebox(0,0)[lt]{\lineheight{1.25}\smash{\begin{tabular}[t]{l}$\mbox{(c)}$\end{tabular}}}}%
  \end{picture}%
\endgroup%

	\label{fig:convergence_vector_linear}
\end{figure}

Moving to the results, Figure \ref{fig:convergence_scalar_linear} shows the $L^{2}$- and $H^{1}$-error norms for the scalar fields $\uP$ and $\lP$ with $k=0$, for both the Natural
 formulation (subfigures (a) and (b)) and the equal-order formulations (subfigures (c)-(f)). In the latter case, the unstabilized, minimally stabilized, and fully stabilized variants are considered. In all cases, the convergence rates are consistent with the expected theoretical error estimate, which is of order $\mathcal{O}(h)$ in the $H^{1}$-norm (see Theorem~2). Also, the error in the $ L^{2} $-norm converges with order $ \mathcal{O}(h^2) $.
No significant differences are observed between the different equal-order variants.
Figure~\ref{fig:convergence_vector_linear} shows the corresponding errors for the vector fields
$\sP$, $\eP$ and $\mP$ in the $ L^{2} $-norm. For the Natural formulation linear convergence is observed for all fields, whereas the equal-order formulation exhibits superlinear behavior, almost reaching $ \mathcal{O}(h^2) $ though the theoretical estimate is $\mathcal{O}(h)$.

\begin{figure}[h!]
	\centering
	\caption{Convergence of scalar fields $\uP$ and $\lP$ with $k=1$. Figures~(a)--(b) show results for $\uP$, and (c)--(d) for $\lP$ in both Natural and equal-order formulations, evaluated in the $L^{2}$ and $H^{1}$ norms.   Three stabilization configurations are compared: Unstabilized, Fully stabilized, and Minimally stabilized ($\alpha$, $\eta$).}
	\def\svgwidth{1.0\columnwidth}
\begingroup%
  \makeatletter%
  \providecommand\color[2][]{%
    \errmessage{(Inkscape) Color is used for the text in Inkscape, but the package 'color.sty' is not loaded}%
    \renewcommand\color[2][]{}%
  }%
  \providecommand\transparent[1]{%
    \errmessage{(Inkscape) Transparency is used (non-zero) for the text in Inkscape, but the package 'transparent.sty' is not loaded}%
    \renewcommand\transparent[1]{}%
  }%
  \providecommand\rotatebox[2]{#2}%
  \newcommand*\fsize{\dimexpr\f@size pt\relax}%
  \newcommand*\lineheight[1]{\fontsize{\fsize}{#1\fsize}\selectfont}%
  \ifx\svgwidth\undefined%
    \setlength{\unitlength}{529.13387269bp}%
    \ifx\svgscale\undefined%
      \relax%
    \else%
      \setlength{\unitlength}{\unitlength * \real{\svgscale}}%
    \fi%
  \else%
    \setlength{\unitlength}{\svgwidth}%
  \fi%
  \global\let\svgwidth\undefined%
  \global\let\svgscale\undefined%
  \makeatother%
  \begin{picture}(1,0.8275543)%
    \lineheight{1}%
    \setlength\tabcolsep{0pt}%
    \put(0.23800499,0.42386128){\color[rgb]{0,0,0}\makebox(0,0)[lt]{\lineheight{1.25}\smash{\begin{tabular}[t]{l}$\mbox{(a)}$\end{tabular}}}}%
    \put(0.73971697,0.42386128){\color[rgb]{0,0,0}\makebox(0,0)[lt]{\lineheight{1.25}\smash{\begin{tabular}[t]{l}$\mbox{(b)}$\end{tabular}}}}%
    \put(0.23862418,0.00320636){\color[rgb]{0,0,0}\makebox(0,0)[lt]{\lineheight{1.25}\smash{\begin{tabular}[t]{l}$\mbox{(c)}$\end{tabular}}}}%
    \put(0.73971697,0.00320636){\color[rgb]{0,0,0}\makebox(0,0)[lt]{\lineheight{1.25}\smash{\begin{tabular}[t]{l}$\mbox{(d)}$\end{tabular}}}}%
    \put(0,0){\includegraphics[width=\unitlength,page=1]{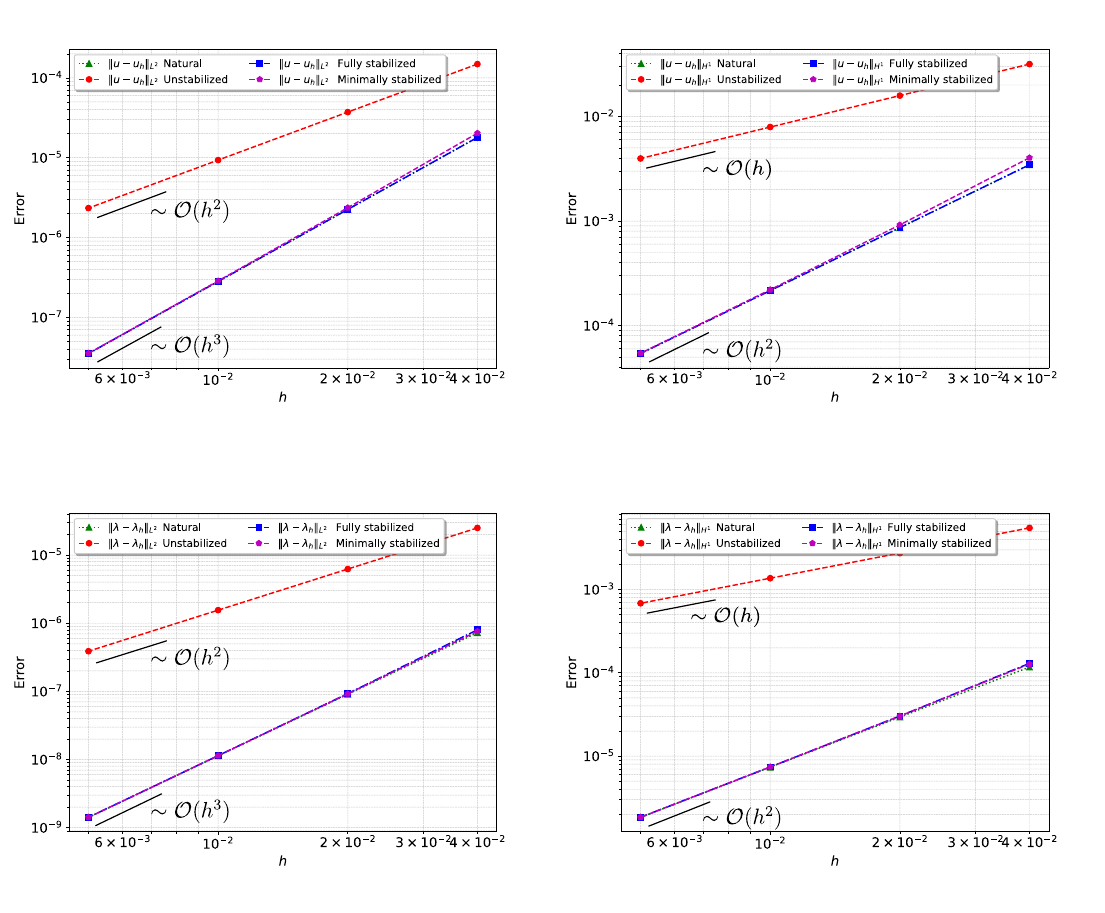}}%
  \end{picture}%
\endgroup%

	\label{fig:convergence_scalar_quadratic}
\end{figure}

\begin{figure}[h!]
	\centering
	\caption{Convergence of vector fields $\sP$, $\eP$, and $\mP$ with $k=1$, evaluated in the $L^{2}$-norm for both Natural and equal-order formulations. Three stabilization configurations are compared: Unstabilized, Fully stabilized, and Minimally stabilized ($\alpha$, $\eta$).}
	\def\svgwidth{1.0\columnwidth}
\begingroup%
  \makeatletter%
  \providecommand\color[2][]{%
    \errmessage{(Inkscape) Color is used for the text in Inkscape, but the package 'color.sty' is not loaded}%
    \renewcommand\color[2][]{}%
  }%
  \providecommand\transparent[1]{%
    \errmessage{(Inkscape) Transparency is used (non-zero) for the text in Inkscape, but the package 'transparent.sty' is not loaded}%
    \renewcommand\transparent[1]{}%
  }%
  \providecommand\rotatebox[2]{#2}%
  \newcommand*\fsize{\dimexpr\f@size pt\relax}%
  \newcommand*\lineheight[1]{\fontsize{\fsize}{#1\fsize}\selectfont}%
  \ifx\svgwidth\undefined%
    \setlength{\unitlength}{529.13387269bp}%
    \ifx\svgscale\undefined%
      \relax%
    \else%
      \setlength{\unitlength}{\unitlength * \real{\svgscale}}%
    \fi%
  \else%
    \setlength{\unitlength}{\svgwidth}%
  \fi%
  \global\let\svgwidth\undefined%
  \global\let\svgscale\undefined%
  \makeatother%
  \begin{picture}(1,0.82205787)%
    \lineheight{1}%
    \setlength\tabcolsep{0pt}%
    \put(0,0){\includegraphics[width=\unitlength,page=1]{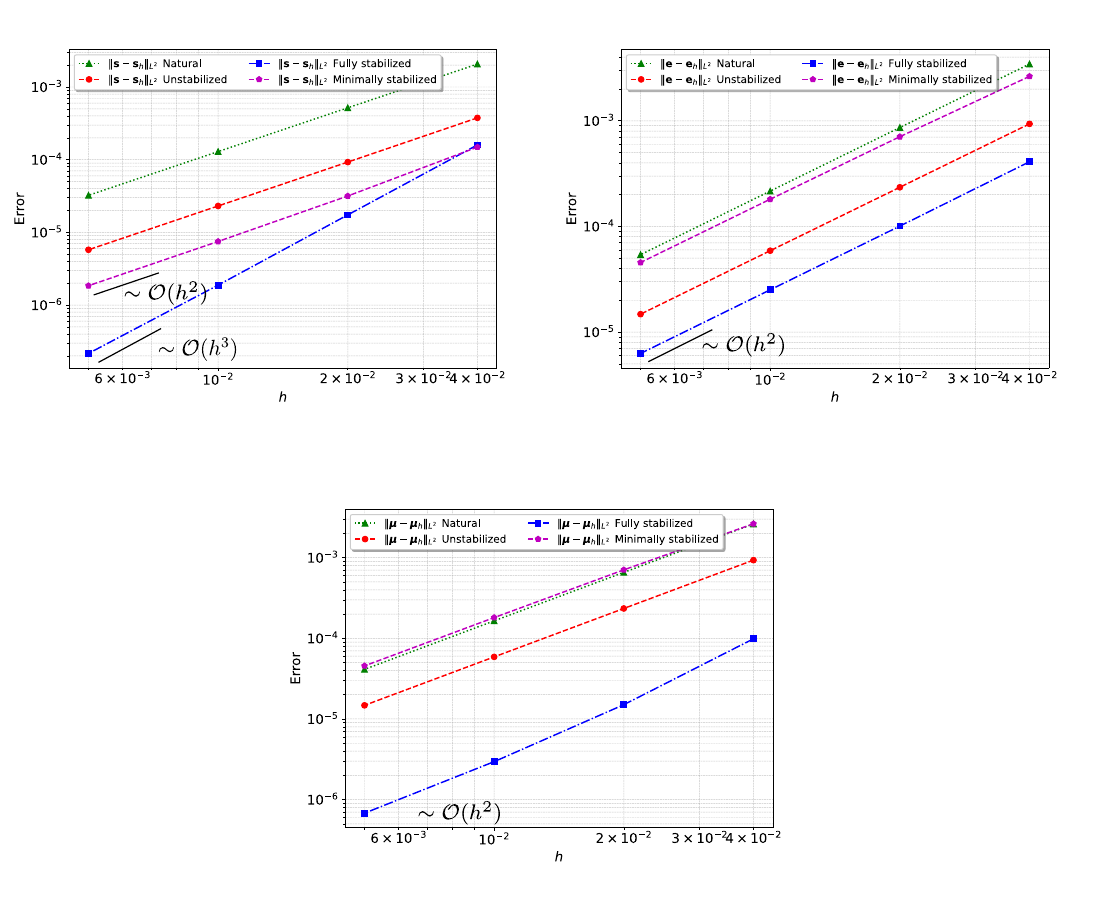}}%
    \put(0.24521587,0.42019699){\color[rgb]{0,0,0}\makebox(0,0)[lt]{\lineheight{1.25}\smash{\begin{tabular}[t]{l}$\mbox{(a)}$\end{tabular}}}}%
    \put(0.74976377,0.42019699){\color[rgb]{0,0,0}\makebox(0,0)[lt]{\lineheight{1.25}\smash{\begin{tabular}[t]{l}$\mbox{(b)}$\end{tabular}}}}%
    \put(0.49391864,0.00320634){\color[rgb]{0,0,0}\makebox(0,0)[lt]{\lineheight{1.25}\smash{\begin{tabular}[t]{l}$\mbox{(c)}$\end{tabular}}}}%
    \put(0,0){\includegraphics[width=\unitlength,page=2]{convergence_vector_quadratic_new.pdf}}%
  \end{picture}%
\endgroup%

	\label{fig:convergence_vector_quadratic}
\end{figure}

\begin{figure}[h!]
	\centering
	\caption{Convergence of scalar fields $\uP$ and $\lP$ with $k=2$. Figures~(a)--(b) show results for $\uP$, and (c)--(d) for $\lP$ in both Natural and equal-order formulations, evaluated in the $L^{2}$ and $H^{1}$ norms.   Three stabilization configurations are compared: Unstabilized, Fully stabilized, and Minimally stabilized ($\alpha$, $\eta$).}
	\def\svgwidth{1.0\columnwidth}
\begingroup%
  \makeatletter%
  \providecommand\color[2][]{%
    \errmessage{(Inkscape) Color is used for the text in Inkscape, but the package 'color.sty' is not loaded}%
    \renewcommand\color[2][]{}%
  }%
  \providecommand\transparent[1]{%
    \errmessage{(Inkscape) Transparency is used (non-zero) for the text in Inkscape, but the package 'transparent.sty' is not loaded}%
    \renewcommand\transparent[1]{}%
  }%
  \providecommand\rotatebox[2]{#2}%
  \newcommand*\fsize{\dimexpr\f@size pt\relax}%
  \newcommand*\lineheight[1]{\fontsize{\fsize}{#1\fsize}\selectfont}%
  \ifx\svgwidth\undefined%
    \setlength{\unitlength}{529.13382943bp}%
    \ifx\svgscale\undefined%
      \relax%
    \else%
      \setlength{\unitlength}{\unitlength * \real{\svgscale}}%
    \fi%
  \else%
    \setlength{\unitlength}{\svgwidth}%
  \fi%
  \global\let\svgwidth\undefined%
  \global\let\svgscale\undefined%
  \makeatother%
  \begin{picture}(1,0.82562277)%
    \lineheight{1}%
    \setlength\tabcolsep{0pt}%
    \put(0.2383659,0.42257361){\color[rgb]{0,0,0}\makebox(0,0)[lt]{\lineheight{1.25}\smash{\begin{tabular}[t]{l}$\mbox{(a)}$\end{tabular}}}}%
    \put(0.74007826,0.42257361){\color[rgb]{0,0,0}\makebox(0,0)[lt]{\lineheight{1.25}\smash{\begin{tabular}[t]{l}$\mbox{(b)}$\end{tabular}}}}%
    \put(0.23898508,0.00320631){\color[rgb]{0,0,0}\makebox(0,0)[lt]{\lineheight{1.25}\smash{\begin{tabular}[t]{l}$\mbox{(c)}$\end{tabular}}}}%
    \put(0.74007826,0.00320631){\color[rgb]{0,0,0}\makebox(0,0)[lt]{\lineheight{1.25}\smash{\begin{tabular}[t]{l}$\mbox{(d)}$\end{tabular}}}}%
    \put(0,0){\includegraphics[width=\unitlength,page=1]{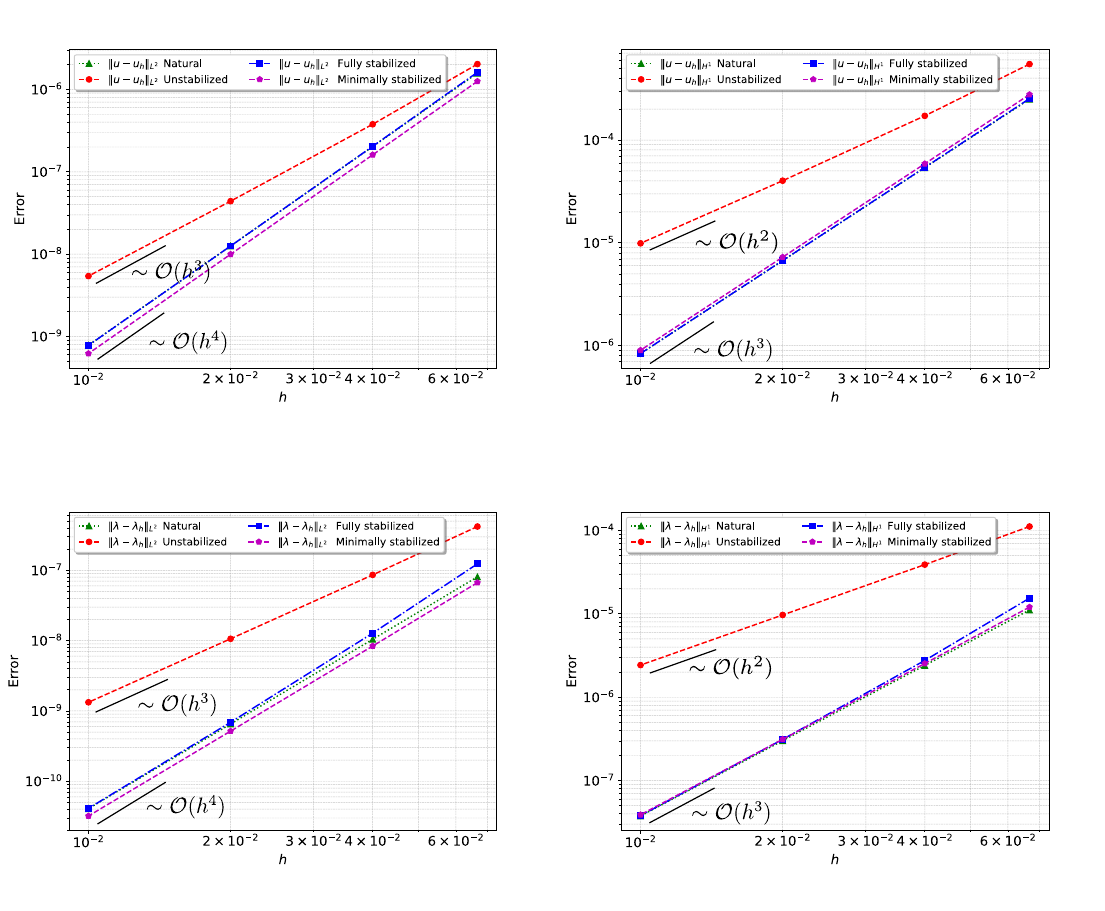}}%
  \end{picture}%
\endgroup%

	\label{fig:convergence_scalar_cubic}
\end{figure}

\begin{figure}[h!]
	\centering
	\caption{Convergence of vector fields $\sP$, $\eP$, and $\mP$ with $k=2$, evaluated in the $L^{2}$-norm for both Natural and equal-order formulations. Three stabilization configurations are compared: Unstabilized, Fully stabilized, and Minimally stabilized ($\alpha$, $\eta$).}
	\def\svgwidth{1.0\columnwidth}
\begingroup%
  \makeatletter%
  \providecommand\color[2][]{%
    \errmessage{(Inkscape) Color is used for the text in Inkscape, but the package 'color.sty' is not loaded}%
    \renewcommand\color[2][]{}%
  }%
  \providecommand\transparent[1]{%
    \errmessage{(Inkscape) Transparency is used (non-zero) for the text in Inkscape, but the package 'transparent.sty' is not loaded}%
    \renewcommand\transparent[1]{}%
  }%
  \providecommand\rotatebox[2]{#2}%
  \newcommand*\fsize{\dimexpr\f@size pt\relax}%
  \newcommand*\lineheight[1]{\fontsize{\fsize}{#1\fsize}\selectfont}%
  \ifx\svgwidth\undefined%
    \setlength{\unitlength}{529.13382943bp}%
    \ifx\svgscale\undefined%
      \relax%
    \else%
      \setlength{\unitlength}{\unitlength * \real{\svgscale}}%
    \fi%
  \else%
    \setlength{\unitlength}{\svgwidth}%
  \fi%
  \global\let\svgwidth\undefined%
  \global\let\svgscale\undefined%
  \makeatother%
  \begin{picture}(1,0.8220581)%
    \lineheight{1}%
    \setlength\tabcolsep{0pt}%
    \put(0,0){\includegraphics[width=\unitlength,page=1]{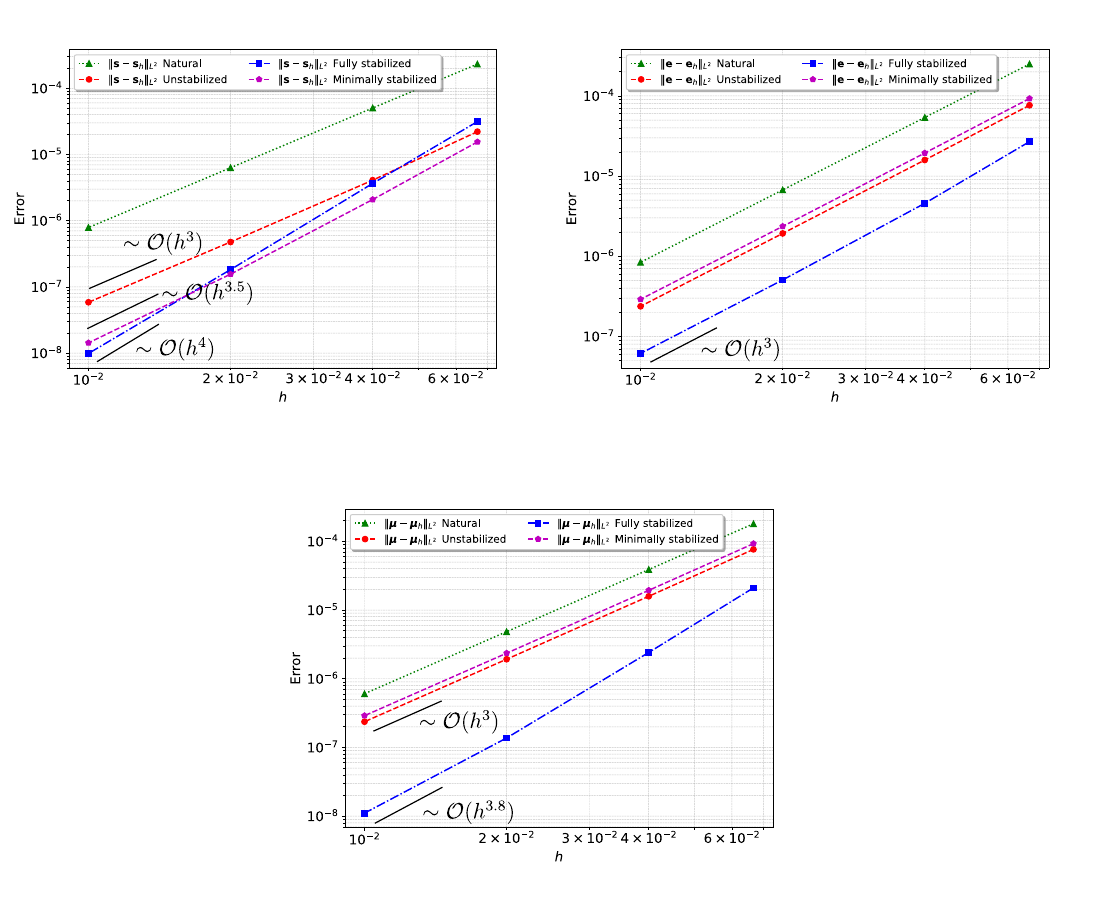}}%
    \put(0.24521595,0.42019721){\color[rgb]{0,0,0}\makebox(0,0)[lt]{\lineheight{1.25}\smash{\begin{tabular}[t]{l}$\mbox{(a)}$\end{tabular}}}}%
    \put(0.74976391,0.42019721){\color[rgb]{0,0,0}\makebox(0,0)[lt]{\lineheight{1.25}\smash{\begin{tabular}[t]{l}$\mbox{(b)}$\end{tabular}}}}%
    \put(0.49391872,0.00320652){\color[rgb]{0,0,0}\makebox(0,0)[lt]{\lineheight{1.25}\smash{\begin{tabular}[t]{l}$\mbox{(c)}$\end{tabular}}}}%
  \end{picture}%
\endgroup%

	\label{fig:convergence_vector_cubic}
\end{figure}

The corresponding numerical results for $\text{P}_2$-elements ($k=1$) and $\text{P}_3$-elements ($k=2$) are shown in Figures~\ref{fig:convergence_scalar_quadratic} and ~\ref{fig:convergence_scalar_cubic}
for the scalar fields and Figures~\ref{fig:convergence_vector_quadratic} and ~\ref{fig:convergence_vector_cubic}
for the vector fields. The Minimally and Fully stabilized 
variants converge at optimal order $\mathcal{O}(h^{k+1})$ for the scalar fields $\uP$ and $\lP$ and for the vector
fields $\eP$ and $\mP$. The convergence rate is better than expected
for the flux $\sP$, attaining $ \mathcal{O}(h^{k+2}) $. Notably, in all cases, the Fully stabilized method exhibited the lowest error values. The lack os stability of the equal-order unstabilized method is evident in the convergence curves of quadratic and cubic elements. They show sub-optimal orders $\mathcal{O}(h^k)$ in the scalar variables $u$ and $\lambda$ in the $H^1$-norm.

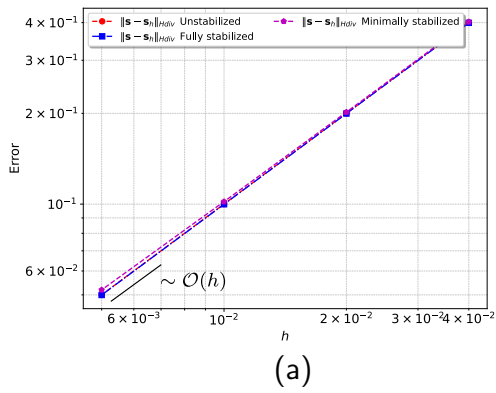
\begin{figure}[h!]
	\centering
	\caption{Convergence of vector fields $\sP$, $\eP$, and $\mP$ with $k=0$ in the equal-order formulation, evaluated in the $\text{Hdiv}$-norm. Figures~(a)--(c) present results for three stabilization configurations: Unstabilized, Fully stabilized, and Minimally stabilized ($\alpha$, $\eta$).
	}
	\def\svgwidth{1.0\columnwidth}
\begingroup%
  \makeatletter%
  \providecommand\color[2][]{%
    \errmessage{(Inkscape) Color is used for the text in Inkscape, but the package 'color.sty' is not loaded}%
    \renewcommand\color[2][]{}%
  }%
  \providecommand\transparent[1]{%
    \errmessage{(Inkscape) Transparency is used (non-zero) for the text in Inkscape, but the package 'transparent.sty' is not loaded}%
    \renewcommand\transparent[1]{}%
  }%
  \providecommand\rotatebox[2]{#2}%
  \newcommand*\fsize{\dimexpr\f@size pt\relax}%
  \newcommand*\lineheight[1]{\fontsize{\fsize}{#1\fsize}\selectfont}%
  \ifx\svgwidth\undefined%
    \setlength{\unitlength}{533.47831666bp}%
    \ifx\svgscale\undefined%
      \relax%
    \else%
      \setlength{\unitlength}{\unitlength * \real{\svgscale}}%
    \fi%
  \else%
    \setlength{\unitlength}{\svgwidth}%
  \fi%
  \global\let\svgwidth\undefined%
  \global\let\svgscale\undefined%
  \makeatother%
  \begin{picture}(1,0.81536342)%
    \lineheight{1}%
    \setlength\tabcolsep{0pt}%
    \put(0.25136252,0.41677511){\color[rgb]{0,0,0}\makebox(0,0)[lt]{\lineheight{1.25}\smash{\begin{tabular}[t]{l}$\mbox{(a)}$\end{tabular}}}}%
    \put(0.75180158,0.41677515){\color[rgb]{0,0,0}\makebox(0,0)[lt]{\lineheight{1.25}\smash{\begin{tabular}[t]{l}$\mbox{(b)}$\end{tabular}}}}%
    \put(0.49803995,0.00318027){\color[rgb]{0,0,0}\makebox(0,0)[lt]{\lineheight{1.25}\smash{\begin{tabular}[t]{l}$\mbox{(c)}$\end{tabular}}}}%
    \put(0,0){\includegraphics[width=\unitlength,page=1]{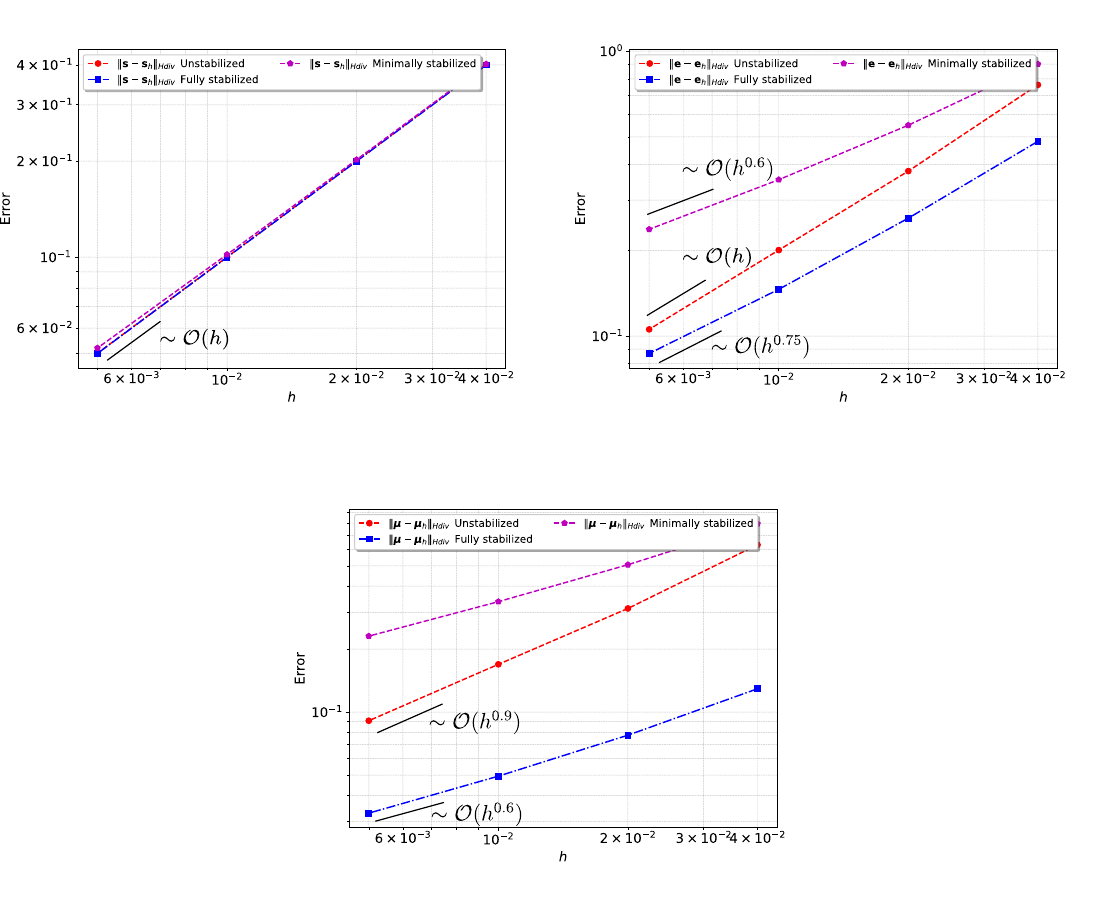}}%
  \end{picture}%
\endgroup%

	\label{fig:convergence_vector_linear_hdiv}
\end{figure}

\begin{figure}[h!]
	\centering
	\caption{Convergence of vector fields $\sP$, $\eP$, and $\mP$ with $k=1$ in the equal-order formulation, evaluated in the $\text{Hdiv}$-norm. Figures~(a)--(c) present results for three stabilization configurations: Unstabilized, Fully stabilized, and Minimally stabilized ($\alpha$, $\eta$).}
	\def\svgwidth{1.0\columnwidth}
\begingroup%
  \makeatletter%
  \providecommand\color[2][]{%
    \errmessage{(Inkscape) Color is used for the text in Inkscape, but the package 'color.sty' is not loaded}%
    \renewcommand\color[2][]{}%
  }%
  \providecommand\transparent[1]{%
    \errmessage{(Inkscape) Transparency is used (non-zero) for the text in Inkscape, but the package 'transparent.sty' is not loaded}%
    \renewcommand\transparent[1]{}%
  }%
  \providecommand\rotatebox[2]{#2}%
  \newcommand*\fsize{\dimexpr\f@size pt\relax}%
  \newcommand*\lineheight[1]{\fontsize{\fsize}{#1\fsize}\selectfont}%
  \ifx\svgwidth\undefined%
    \setlength{\unitlength}{529.13387269bp}%
    \ifx\svgscale\undefined%
      \relax%
    \else%
      \setlength{\unitlength}{\unitlength * \real{\svgscale}}%
    \fi%
  \else%
    \setlength{\unitlength}{\svgwidth}%
  \fi%
  \global\let\svgwidth\undefined%
  \global\let\svgscale\undefined%
  \makeatother%
  \begin{picture}(1,0.82130264)%
    \lineheight{1}%
    \setlength\tabcolsep{0pt}%
    \put(0.23833753,0.41969356){\color[rgb]{0,0,0}\makebox(0,0)[lt]{\lineheight{1.25}\smash{\begin{tabular}[t]{l}$\mbox{(a)}$\end{tabular}}}}%
    \put(0.74005085,0.41969356){\color[rgb]{0,0,0}\makebox(0,0)[lt]{\lineheight{1.25}\smash{\begin{tabular}[t]{l}$\mbox{(b)}$\end{tabular}}}}%
    \put(0.48987512,0.00320637){\color[rgb]{0,0,0}\makebox(0,0)[lt]{\lineheight{1.25}\smash{\begin{tabular}[t]{l}$\mbox{(c)}$\end{tabular}}}}%
    \put(0,0){\includegraphics[width=\unitlength,page=1]{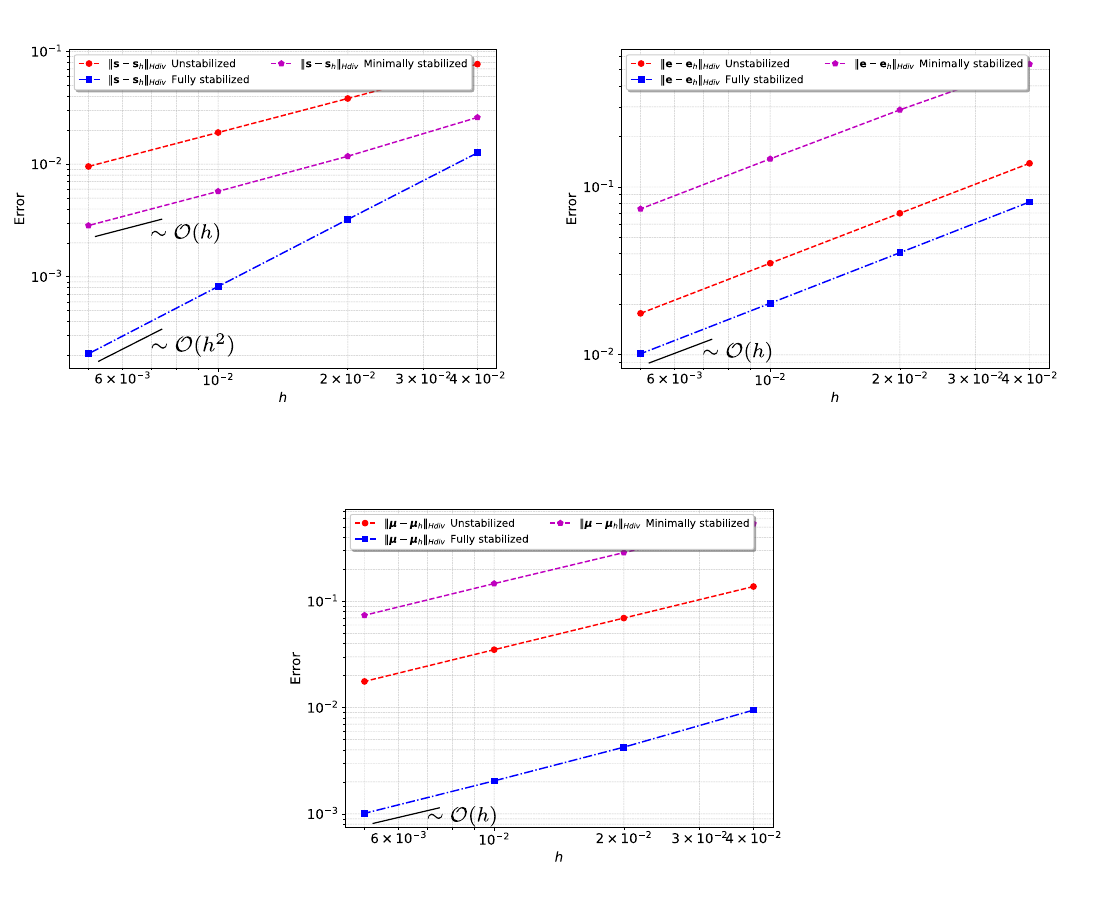}}%
  \end{picture}%
\endgroup%

	\label{fig:convergence_vector_quadratic_hdiv}
\end{figure}

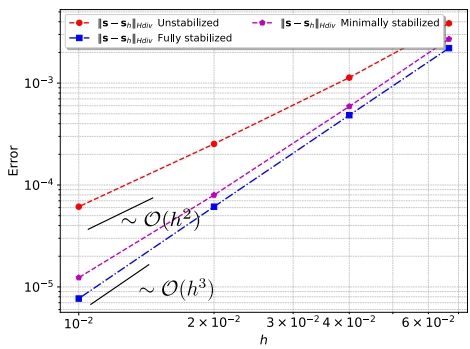
\begin{figure}[h!]
	\centering
	\caption{Convergence of vector fields $\sP$, $\eP$, and $\mP$ with $k=2$ in the equal-order formulation, evaluated in the $\text{Hdiv}$-norm. Figures~(a)--(c) present results for three stabilization configurations: Unstabilized, Fully stabilized, and Minimally stabilized ($\alpha$, $\eta$).}
	\def\svgwidth{1.0\columnwidth}
\begingroup%
  \makeatletter%
  \providecommand\color[2][]{%
    \errmessage{(Inkscape) Color is used for the text in Inkscape, but the package 'color.sty' is not loaded}%
    \renewcommand\color[2][]{}%
  }%
  \providecommand\transparent[1]{%
    \errmessage{(Inkscape) Transparency is used (non-zero) for the text in Inkscape, but the package 'transparent.sty' is not loaded}%
    \renewcommand\transparent[1]{}%
  }%
  \providecommand\rotatebox[2]{#2}%
  \newcommand*\fsize{\dimexpr\f@size pt\relax}%
  \newcommand*\lineheight[1]{\fontsize{\fsize}{#1\fsize}\selectfont}%
  \ifx\svgwidth\undefined%
    \setlength{\unitlength}{529.13387269bp}%
    \ifx\svgscale\undefined%
      \relax%
    \else%
      \setlength{\unitlength}{\unitlength * \real{\svgscale}}%
    \fi%
  \else%
    \setlength{\unitlength}{\svgwidth}%
  \fi%
  \global\let\svgwidth\undefined%
  \global\let\svgscale\undefined%
  \makeatother%
  \begin{picture}(1,0.82205787)%
    \lineheight{1}%
    \setlength\tabcolsep{0pt}%
    \put(0.24521589,0.42019703){\color[rgb]{0,0,0}\makebox(0,0)[lt]{\lineheight{1.25}\smash{\begin{tabular}[t]{l}$\mbox{(a)}$\end{tabular}}}}%
    \put(0.74976388,0.42019703){\color[rgb]{0,0,0}\makebox(0,0)[lt]{\lineheight{1.25}\smash{\begin{tabular}[t]{l}$\mbox{(b)}$\end{tabular}}}}%
    \put(0.49391866,0.0032063){\color[rgb]{0,0,0}\makebox(0,0)[lt]{\lineheight{1.25}\smash{\begin{tabular}[t]{l}$\mbox{(c)}$\end{tabular}}}}%
    \put(0,0){\includegraphics[width=\unitlength,page=1]{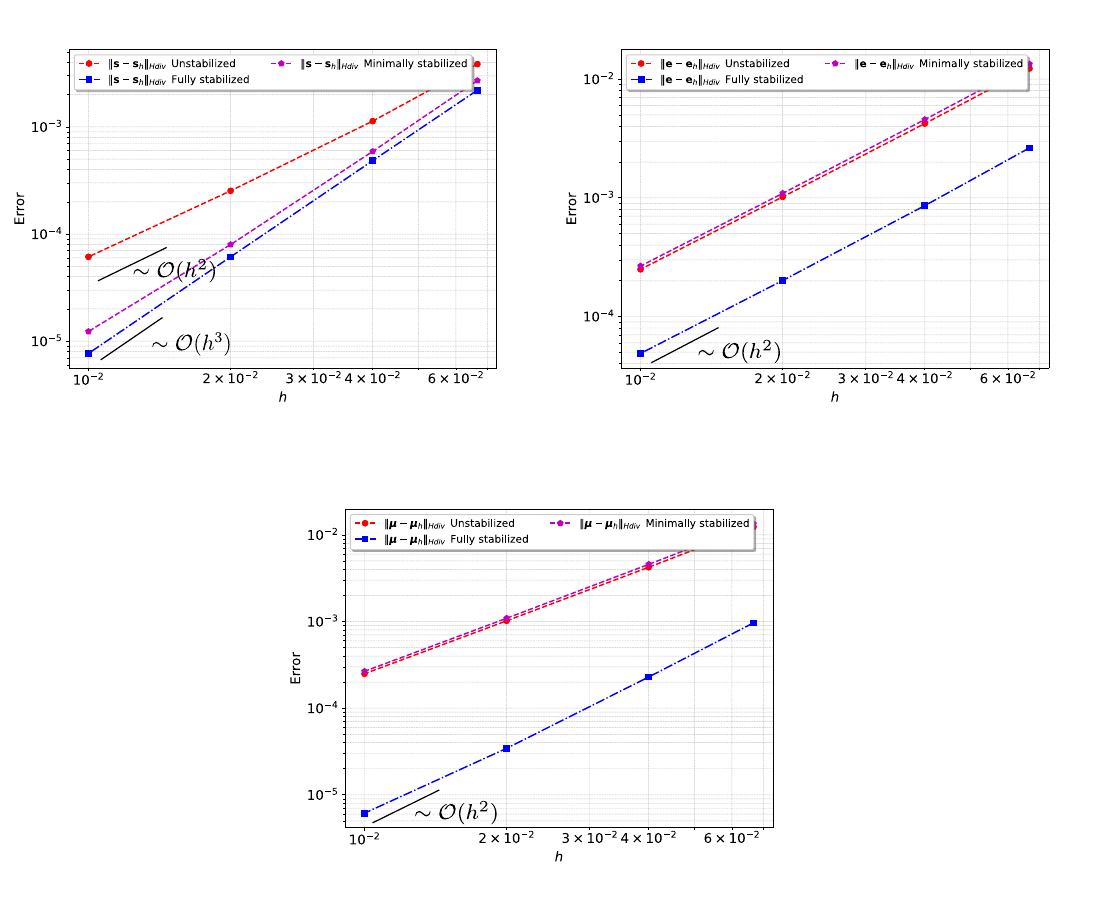}}%
  \end{picture}%
\endgroup%

	\label{fig:convergence_vector_cubic_hdiv}
\end{figure}

Although we do not have theoretical error estimates in the $ \text{H(div)} $-norm, 
convergence rates are also reported. Figures \ref{fig:convergence_vector_linear_hdiv}, 
\ref{fig:convergence_vector_quadratic_hdiv} and \ref{fig:convergence_vector_cubic_hdiv}
show the results for $k=0$, $k=1$ and $k=2$, respectively.
Note that the Fully stabilized formulation yielded the most favorable results in all three cases. For $ \text{P}_1 $-elements, no significant differences were observed in the convergence behavior of the flux $ \sP_{h} $ across the Unstabilized, Minimally stabilized, and Fully stabilized formulations, with the corresponding error curves largely overlapping. In contrast, for $ \eP $ and $ \mP $, the Fully stabilized formulation exhibited sub-linear convergence, unlike the Unstabilized case. However, it achieved lower error magnitudes. For $ \text{P}_2 $-elements, the Fully stabilized formulation exhibited the best performance. In particular, it was the only one to attain second-order convergence for the flux $ \sP $, and exhibited the smallest errors for $ \eP_{h} $ and $ \mP_{h} $. A similar trend was observed for $ \text{P}_3 $-elements. 

It is worth mentioning that additional numerical experiments were carried out to evaluate the effect of each individual stabilization parameter in isolation and compared to the Unstabilized, Minimally stabilized, and Fully stabilized cases. The results are not reported here for the sake of brevity, however, they displayed only minor variations in magnitude and exhibited convergence rates consistent with those of the main reference cases. Furthermore, numerical experiments confirmed that the convergence behavior is rather independent of the reaction coefficient $ \zeta $, as similar results were obtained when varying $\zeta$ over several orders of magnitude.

\begin{figure}[h!]
	\centering
	\caption{Figures~(a)--(d) show the difference $\uP_{h} - \uP$ for the Natural formulation, and for the equal-order Unstabilized, Fully stabilized, and Minimally stabilized formulations, respectively. }
	\def\svgwidth{1.0\columnwidth}
	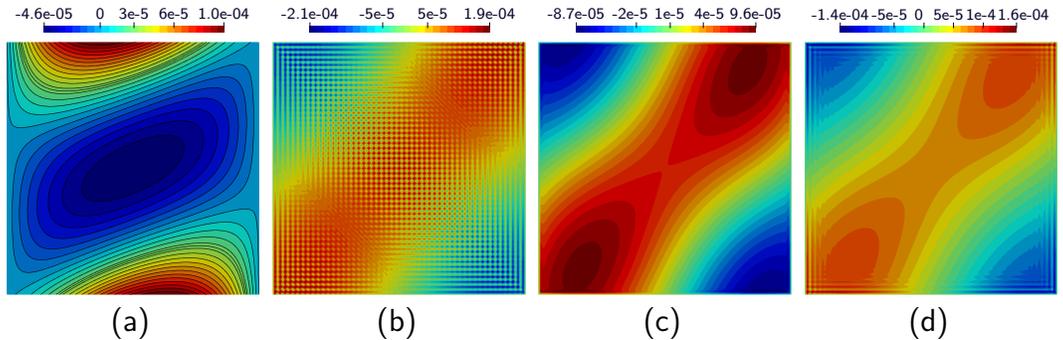
	\label{fig:checkerboard_pattern}
\end{figure}

In what follows, we illustrate the stability properties of the approximation for the Natural and 
the equal-order formulations by evaluating the difference $\uP_{h} - \uP$ on a sufficiently fine mesh. Figure~\ref{fig:checkerboard_pattern} presents the corresponding results. We observe that no spurious modes 
are present in either the Natural formulation or the Equal-order Fully stabilized method. 
In the Minimally stabilized method, a subtle checkerboard pattern appears near the boundaries. In contrast, in the Unstabilized method, spurious modes are manifested throughout the entire 
domain as expected. However, it is worth mentioning that this pattern becomes evident 
only when we plot the pointwise values of the error $\uP_{h} - \uP$ and the norm of 
this error tends to zero (with suboptimal order) when the mesh is refined.

As a final comment, we emphasize that the manufactured fields $ \eP $ and $ \sP $ were constructed to be physically admissible in the sense that they do not violate the second law of thermodynamics 
\cite{Reddy2013,Gurtin2013}. Accordingly, it holds that
\begin{equation}
\sP(\mathbf{x}) \cdot \eP(\mathbf{x}) \leq 0, \quad \forall \mathbf{x} \in \Omega
\label{eq:duhem}
\end{equation}
Even though the discrete approximations $ \sP_{h} $ and $ \eP_{h} $ converge to the manufactured fields in norm, this does not guarantee that $\sP_h \cdot \eP_h$ is negative everywhere. Therefore, for all problems solved in this work we have explicitly computed the sign of $\sP_h \cdot \eP_h$ at the finite element nodes. Remarkably, no violations of (\ref{eq:duhem}) were detected. 


\cleardoublepage
\newpage
\subsection{Test case 2: a non-convex domain}\label{sec:test_case2}

\begin{figure}[h!]
	\centering
	\caption{Illustration of the three non-convex geometries considered.}
	\fontsize{11pt}{10pt}\selectfont 
	\def\svgwidth{1.0\columnwidth}
	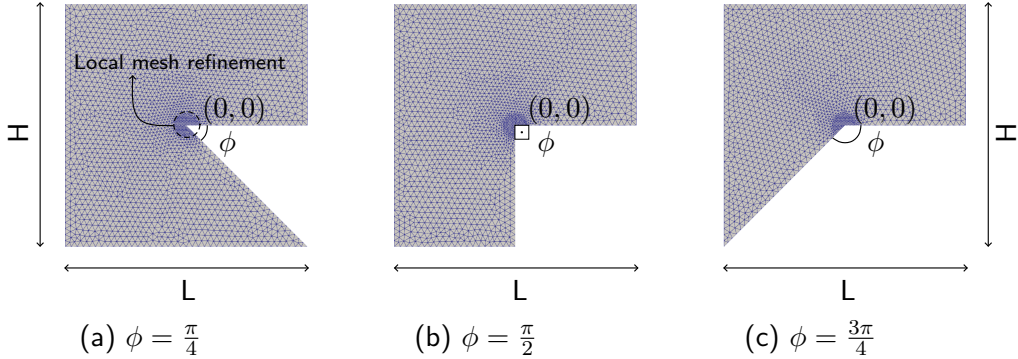
	\label{fig:geometries_non_convex}
\end{figure}

\begin{figure}[h!]
	\centering
	\caption{Figure (a) presents the functional relationship between the norm of the flux and the gradient $\|\sP\|~\mbox{vs.}~\|\eP\|$, and $\|\stP\|~\mbox{vs.}~\|\etP\|$, respectively. Figure (b) to (f) show $\uP$, $\lP$, $\eP$, $\sP$, and $\mP$.}
	\def\svgwidth{1.0\columnwidth}
	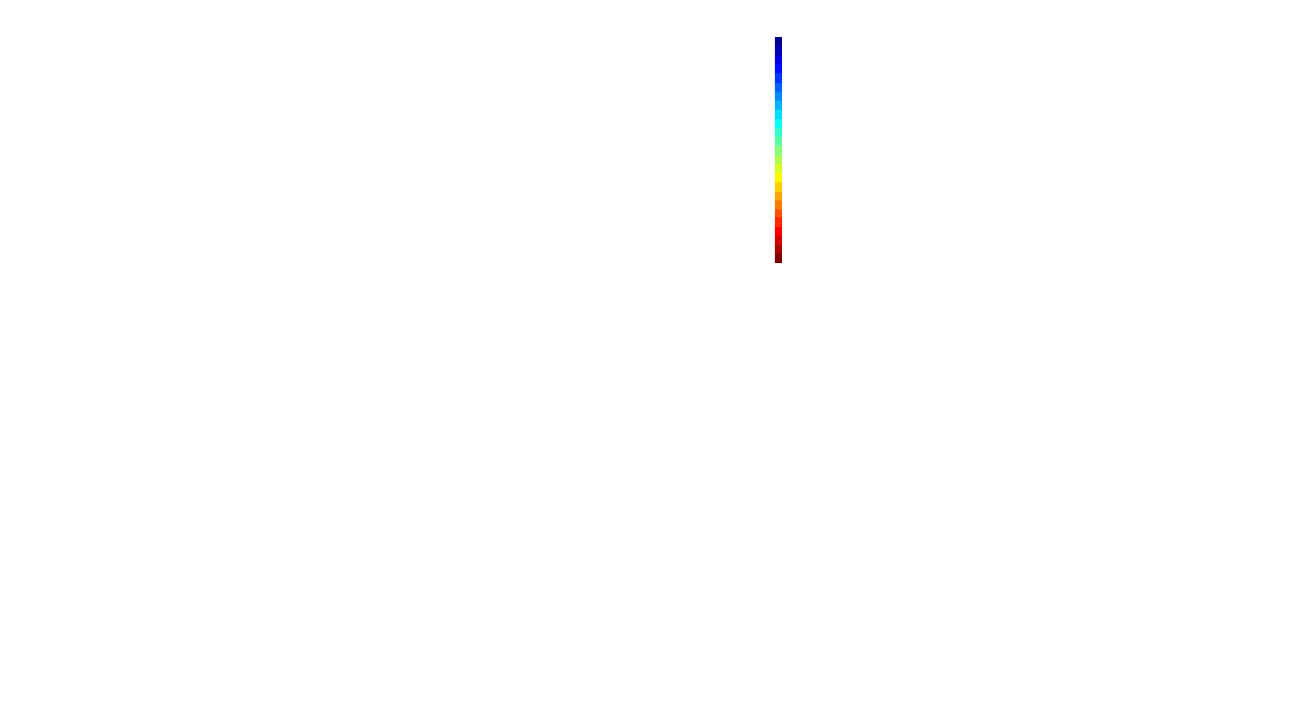
	\label{fig:all_fields_ref_non_convex}
\end{figure}

In this case we consider the classical problem with a re-entrant corner, whose geometry is depicted in 
Figure \ref{fig:geometries_non_convex} for three different external angles, namely, $\phi =\pi/4,~\pi/2$ and $3\pi/4$. The corner is located at the origin $(0, 0)$. 
We manufacture a solution by taking the scalar field $u$ as being:
\begin{equation*}
	\begin{aligned}
	&\uP: \Omega \rightarrow \mathbb{R}, \quad \uP(x,y) = r^{\nu}\,\sin(\nu\,\theta) \\
	&\mbox{where}~r = \sqrt{\left(x^{2}+y^{2}\right)}, \quad \theta = \arctan{\left(\frac{y}{x}\right)}, \quad  \nu = \frac{\pi}{\psi},~\mbox{with}~\psi = 2\pi-\phi
	\end{aligned}
\end{equation*}
The corresponding gradient field is then given by
\begin{equation*}
	\begin{aligned}
	\eP: \Omega \rightarrow \mathbb{R}^{2}, \quad \eP(x,y) = \nu\,r^{\nu-2}\left[x\,\sin(\nu\,\theta)-y\,\cos(\nu\,\theta),y\,\sin(\nu\,\theta)+x\,\cos(\nu\,\theta)\right]^{\top}
	\end{aligned}
\end{equation*}
In this setting, the expressions for $\lP$ and $\etP$ are taken to be the same as those used in the preceding numerical experiment. For $\sP$ we simply consider the linear relation $\sP = -\eP$.
Also, we assume pure Dirichlet boundary conditions for $\uP$ and $\lP$ based on the manufactured solution.
Figure~\ref{fig:all_fields_ref_non_convex} illustrates the manufactured fields. Again, we also show the behavior of the norm of $\sP$ as a function of the norm of  $\eP$ and their counterparts $ \stP$ and $\etP $.
It is worth noticing that the cubic term considered in the previous experiment has been suppressed so as to not exacerbate the singular behavior of the solution.
For all these geometries, it can be shown that $\uP$ possesses fewer than $\nu+1$ derivatives in $L^{2}$, implying that $\uP \in H^{1}(\Omega)$, but $\uP \notin H^{2}(\Omega)$. Consequently, no improvement in the convergence rate is expected from polynomial refinement and therefore we only use $\text{P}_1$-elements.
In order to better capture the solution, a local mesh refinement near the singularity is employed, as illustrated in Figure~\ref{fig:geometries_non_convex}(a).
We assess the convergence behavior with respect to different corner angles. Specifically, for the values $ \phi \in \left\{ \frac{\pi}{4}, \frac{\pi}{2}, \frac{3\pi}{4} \right\} $, the corresponding values of $ \psi $ are $ \left\{ \frac{7\pi}{4}, \frac{3\pi}{2}, \frac{5\pi}{4} \right\} $. Consequently, the scalar field $ \uP $ belong to the spaces $ H^{11/7 - \varepsilon}(\Omega) $, $ H^{5/3 - \varepsilon}(\Omega) $, and $ H^{9/5 - \varepsilon}(\Omega) $, for all $ \varepsilon > 0 $. Error estimates for finite element schemes in domains with corners have been extensively studied in the literature, for instance by \cite{Babuska1972, Babuska1970}. As a result, the convergence rate of $ \uP $ in both $ L^{2} $-norm  and $ H^{1} $-norms are expected to be restricted by the singular behavior at the reentrant corner and one may anticipate convergence rates limited by $\mathcal{O}(h^{1+4/7-\varepsilon})$, $ \mathcal{O}(h^{1+2/3-\varepsilon}) $ and $ \mathcal{O}(h^{1+4/5-\varepsilon}) $ for all $ \varepsilon > 0 $ in $L^{2}$-norm, and $\mathcal{O}(h^{4/7-\varepsilon})$, $ \mathcal{O}(h^{2/3-\varepsilon}) $ and $ \mathcal{O}(h^{4/5-\varepsilon}) $ for all $ \varepsilon > 0 $ in $H^{1}$-norm. 
Table~\ref{tab:scalar_errors} summarizes the convergence rates for the scalar fields in the $ H^{1} $-norms, considering the Natural formulation and the Equal-order formulations, across the three corner angles. The suboptimal rates seen for the Unstabilized method are highlighted.

Similarly, Table~\ref{tab:vector_errors} reports the convergence rates of the vector fields in the $ L^{2} $-norm.
Focusing on the scalar field $ \uP $ in the Natural formulation, we observe that the convergence in the $ H^1 $-norm closely approaches the theoretical upper bounds. The field $ \lP $ exhibits slightly better convergence than $ \uP $, particularly in the $ H^1 $-norm.
\begin{table}[h!]
	\centering
	\caption{Estimated convergence rates for the scalar fields $ \uP_h $ and $ \lP_h $ in the $ L^2 $- and $ H^1 $-norms using different stabilization schemes and reentrant corner angles $ \psi \in \left\{ \frac{7\pi}{4}, \frac{3\pi}{2}, \frac{5\pi}{4} \right\} $.}
	\begin{tabular}{l l l l l}
		\toprule
		Formulation & Norms & $\psi = \frac{7\pi}{4}$ & $\psi = \frac{3\pi}{2}$ & $\psi = \frac{5\pi}{4}$ \\ [2pt]
		\hline
		Interpolation                   &  & 0.57 & 0.66 & 0.80 \\  
		\hline
		\multirow{2}{*}{Natural}		& $\|u-\uP_{h}\|_{H^{1}}$ & 0.57 & 0.72  & 0.81 \\
		                                 & $\|\lP-\lP_{h}\|_{H^{1}}$ & 0.99  & 0.99 & 0.99  \\
		\hline
		\multirow{2}{*}{Unstabilized}    & $\|u-\uP_{h}\|_{H^{1}}$ & \textbf{0.46} & \textbf{0.59} & \textbf{0.64} \\
		                                 & $\|\lP-\lP_{h}\|_{H^{1}}$ & \textbf{0.46} & \textbf{0.30} & 1.24 \\
		\hline
		\multirow{2}{*}{Minimally stabilized} & $\|u-\uP_{h}\|_{H^{1}}$ & 0.57 & 0.73 & 0.81 \\
		& $\|\lP-\lP_{h}\|_{H^{1}}$ & 0.60 & 0.78 & 0.98 \\
		\hline
		\multirow{2}{*}{Fully stabilized} & $\|u-\uP_{h}\|_{H^{1}}$ & 0.57 & 0.72 & 0.81 \\
		& $\|\lP-\lP_{h}\|_{H^{1}}$ & 0.60 & 0.82 & 0.99 \\
		\bottomrule
	\end{tabular}
	\label{tab:scalar_errors}
\end{table}

When comparing the same fields under the equal-order formulation, the convergence rates remain similar in magnitude to those obtained with the Natural formulation, with no substantial improvement. Among the equal-order variants, we highlight the enhanced convergence of $ \uP $ and $ \lP $ in the $ H^1 $-norm achieved by the Minimally and Fully stabilized formulations, compared to the Unstabilized case. However, no further improvements are observed in the $ L^2 $-norm for either field when stabilization is applied.
For the vector fields $ \sP$, $ \eP $, and $ \mP $, we observe that the convergence rates in the $ L^2 $-norm remain nearly identical across all formulations and geometries. Improvements in convergence are primarily associated with increasing the interior angle at the re-entrant corner, with higher accuracy achieved when transitioning from acute to obtuse geometries. 
These results confirm that the convergence behavior is primarily dictated by the corner-induced singularity. While stabilization enhances robustness -- especially in the $H^1$-norm for scalar fields -- it does not overcome the inherent regularity limitations imposed by the geometry. 
\begin{table}[h!]
	\centering
	\caption{Estimated convergence rates in the $ L^2 $-norm for the vector fields $ \sP_h $, $ \eP_h $, and $ \mP_h $ using various stabilization schemes and different reentrant corner angles $ \psi \in \left\{ \frac{7\pi}{4}, \frac{3\pi}{2}, \frac{5\pi}{4} \right\} $.}
	\begin{tabular}{l l c c c}
		\toprule
		Formulation & Norms & $\psi = \frac{7\pi}{4}$ & $\psi = \frac{3\pi}{2}$ & $\psi = \frac{5\pi}{4}$ \\ [2pt]
		\hline
		Interpolation                   &  & 0.57 & 0.66 & 0.80 \\
		\hline  
		\multirow{3}{*}{Natural}    & $\|\mathbf{s}-\sP_{h}\|_{L^{2}}$ & 0.57 & 0.73 & 0.83  \\
		                                 & $\|\mathbf{e}-\eP_{h}\|_{L^{2}}$ & 0.57 & 0.72 & 0.81 \\
		                                 & $\|\mP-\mP_{h}\|_{L^{2}}$ & 0.58 & 0.73 & 0.87 \\
		\hline
		\multirow{3}{*}{Unstabilized}    & $\|\mathbf{s}-\sP_{h}\|_{L^{2}}$ & 0.55 & 0.67 & 0.86 \\
		                                 & $\|\mathbf{e}-\eP_{h}\|_{L^{2}}$ & 0.56 & 0.72 & 0.84 \\
		                                 & $\|\mP-\mP_{h}\|_{L^{2}}$ & 0.64 & 0.69 & 0.97 \\
		\hline
				\multirow{3}{*}{Minimally stabilized} & $\|\mathbf{s}-\sP_{h}\|_{L^{2}}$ & 0.55 & 0.73 & 0.82 \\
		& $\|\mathbf{e}-\eP_{h}\|_{L^{2}}$ & 0.55 & 0.72 & 0.82 \\
		& $\|\mP-\mP_{h}\|_{L^{2}}$ & 0.58 & 0.69 & 0.83 \\
		\hline
		\multirow{3}{*}{Fully stabilized} & $\|\mathbf{s}-\sP_{h}\|_{L^{2}}$ & 0.56 & 0.72 & 0.84 \\
		& $\|\mathbf{e}-\eP_{h}\|_{L^{2}}$ & 0.55 & 0.73 & 0.82 \\
		& $\|\mP-\mP_{h}\|_{L^{2}}$ & 0.55 & 0.73 & 0.89 \\
		\bottomrule
	\end{tabular}
	\label{tab:vector_errors}
\end{table}

\newpage
\subsection{Test case 3: data assignation}

In the previous numerical examples, the exact values of $ (\etP, \stP) $ were available at every integration point in $\Omega$, such that we could assess the convergence
properties of the finite element formulation. 
Now, as a final numerical experiment, we assess the consistency of the discrete 
formulation with respect to the number of available data pairs $ (\etP, \stP) $, aiming 
to mimic a more realistic scenario for data-driven applications. 

To this end, we consider a structured cartesian grid of the computational domain into $N_{\mbox{\tiny samples}} = N_d \times N_d$ non-overlapping quadrilateral subdomains $ \{\Omega_i\}_{i=1}^{N_{\mbox{\tiny samples}}} $ where the solution will be sampled. For each subdomain $ \Omega_i $, let $\mathbf{x}_{i}^{c} \in \Omega$ denotes its centroid, at which a unique data pair $(\etP_i, \stP_i)$ is associated. For this problem we consider the manufactured solution: 
\begin{equation}
	\uP: \Omega \rightarrow \mathbb{R}, \quad \uP(x,y) = \sin(\pi x)\sin(\pi y),
\end{equation}
whose gradient is given by:
\begin{equation}
	\eP: \Omega \rightarrow \mathbb{R}^{2}, \quad \eP(x,y) = \pi \, \left[\cos(\pi x)\sin(\pi y), \sin(\pi x)\cos(\pi y) \right]^{\top}.
\end{equation}
For the sake of simplicity, we consider $ \sP = -\eP $ and further assume $\lP  = 0$, $ \mP = \left[0,0\right]^{\top}$,
such that, the data fields $(\etP, \stP)$ coincide with $(\eP, \sP)$ in the exact setting. 
Additionally, for the source term in Eq.~(\ref{eq:dif4}) we have: 
$$
q: \Omega \rightarrow \mathbb{R}, \quad q(x,y) = \left(1 + 2\pi^2 \right)\sin(\pi x)\sin(\pi y).
$$

The values of $\{(\etP_i, \stP_i)\}_{i=1}^{N_{\mbox{\tiny samples}}}$ are thus obtained by evaluating the manufactured fields $(\eP, \sP)$ at $\mathbf{x}_{i}^{c}$. To assign a data pair $ (\etP_i, \stP_i) $ from the dataset to each element of the mesh, we use a geometric proximity criterion, i.e., denoting by $ \mathbf{x}_K^{c} $ the centroid of mesh 
element $ K \in \mathcal{T}_h $, the data point assigned to element $K$ is given by
\begin{equation}
	(\etP_h, \stP_h)(\mathbf{x}_{K}^{c}) := (\etP_i, \stP_i) \quad \text{where} \quad i = \arg\min_{j} \| \mathbf{x}_K^{c} - \mathbf{x}_{j}^{c} \|_2
\end{equation}
where $\| \cdot \|_2$ stands for the Euclidean norm. Notice that, in this way, the data fields $(\etP, \stP)$
have an elementwise constant representation over $\Omega$ in which all the elements of
$\mathcal{T}_h$ that are close to a sampled point share the same value. 
This does not introduce any additional
approximation error only when the vector fields are in $[\DG_0]^2$.
Finally, in the numerical experiments, we consider only $ k=0 $ for both formulations 
and assume pure Dirichlet boundary conditions.

\begin{figure}[h!]
	\centering
	\caption{Contours of pressure and flux obtained using the Natural formulation, with an increasing number of data samples distributed over the domain. The last row also shows the corresponding flux field computed using the equal-order formulation. The pressure field for this case is omitted, as it is very
	similar to that obtained with the Natural formulation.}
	\fontsize{8pt}{10pt}\selectfont 
	\def\svgwidth{1.0\columnwidth}
	\input{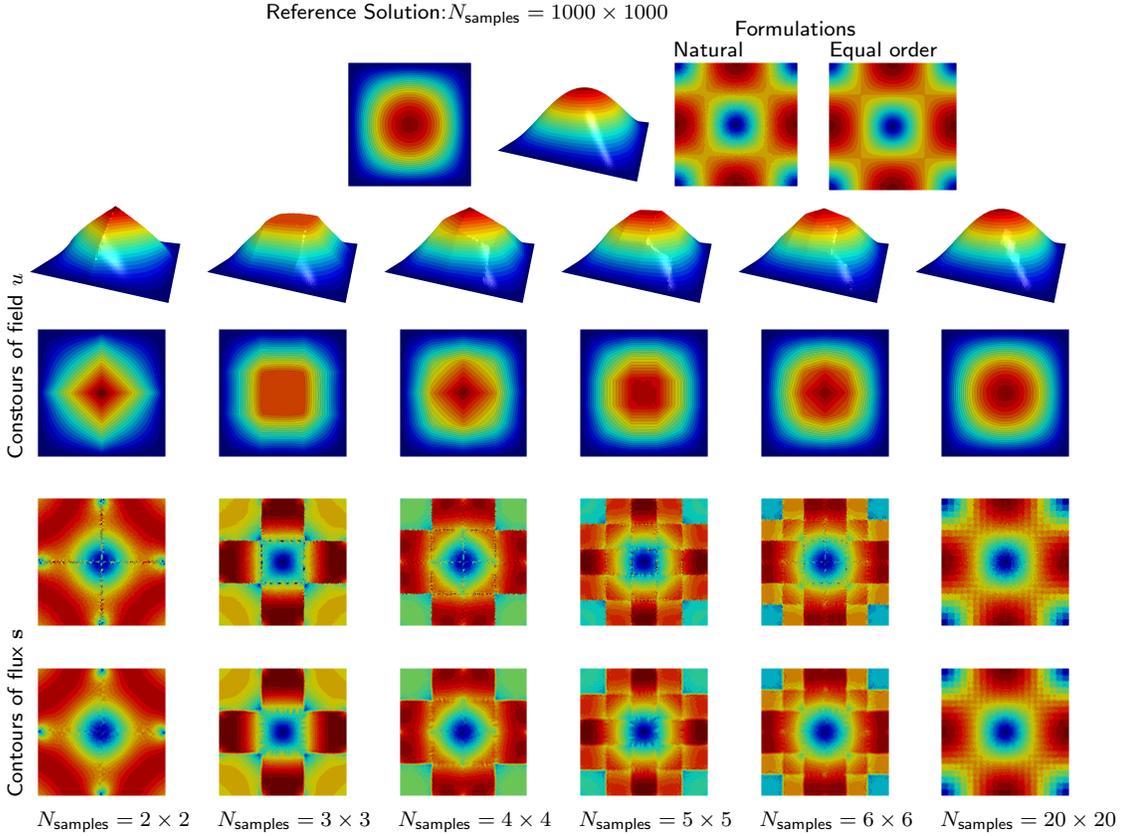}
	\label{fig:compound_Nsamples}
\end{figure}

To illustrate the solution obtained with this approach, figure \ref{fig:compound_Nsamples} shows
contour plots of pressure and flux for the Natural formulation 
varying the number of subdomains from $2\times 2$ samples up to $20 \times 20$. 
In the last row of figures, we also show contours of flux for the equal-order
formulation, the pressure for this case being very similar to the former. 
The results can be compared to the most refined partition, which encompasses $1000 \times 1000$
samples and is indistinguishable from the exact solution. 
The coarsest samplings produce pressure fields that exhibit a pyramidal shape with
as many faces as number of samples considered. The $20 \times 20$ solution 
already shows a high degree of similarity to the exact one.

\begin{figure}[h!]
	\centering
	\caption{Mesh convergence results for the scalar fields $\uP$ and $\lP$ in the $L^2$ and $H^1$ norms, and for the vector fields $\eP$, $\sP$, and $\mP$ in the $L^2$ norm, using the Natural
	 formulation across various numbers of subdomains $N_{\mbox{\tiny samples}}$.}
	\def\svgwidth{1.0\columnwidth}
\begingroup%
  \makeatletter%
  \providecommand\color[2][]{%
    \errmessage{(Inkscape) Color is used for the text in Inkscape, but the package 'color.sty' is not loaded}%
    \renewcommand\color[2][]{}%
  }%
  \providecommand\transparent[1]{%
    \errmessage{(Inkscape) Transparency is used (non-zero) for the text in Inkscape, but the package 'transparent.sty' is not loaded}%
    \renewcommand\transparent[1]{}%
  }%
  \providecommand\rotatebox[2]{#2}%
  \newcommand*\fsize{\dimexpr\f@size pt\relax}%
  \newcommand*\lineheight[1]{\fontsize{\fsize}{#1\fsize}\selectfont}%
  \ifx\svgwidth\undefined%
    \setlength{\unitlength}{533.47831666bp}%
    \ifx\svgscale\undefined%
      \relax%
    \else%
      \setlength{\unitlength}{\unitlength * \real{\svgscale}}%
    \fi%
  \else%
    \setlength{\unitlength}{\svgwidth}%
  \fi%
  \global\let\svgwidth\undefined%
  \global\let\svgscale\undefined%
  \makeatother%
  \begin{picture}(1,1.57359276)%
    \lineheight{1}%
    \setlength\tabcolsep{0pt}%
    \put(0.23537178,0.37964164){\color[rgb]{0,0,0}\makebox(0,0)[lt]{\lineheight{1.25}\smash{\begin{tabular}[t]{l}$\mbox{(e)}$\end{tabular}}}}%
    \put(0.2353956,1.18136988){\color[rgb]{0,0,0}\makebox(0,0)[lt]{\lineheight{1.25}\smash{\begin{tabular}[t]{l}$\mbox{(a)}$\end{tabular}}}}%
    \put(0.73518054,1.18136988){\color[rgb]{0,0,0}\makebox(0,0)[lt]{\lineheight{1.25}\smash{\begin{tabular}[t]{l}$\mbox{(b)}$\end{tabular}}}}%
    \put(0.23600975,0.78050572){\color[rgb]{0,0,0}\makebox(0,0)[lt]{\lineheight{1.25}\smash{\begin{tabular}[t]{l}$\mbox{(c)}$\end{tabular}}}}%
    \put(0.73518054,0.78050572){\color[rgb]{0,0,0}\makebox(0,0)[lt]{\lineheight{1.25}\smash{\begin{tabular}[t]{l}$\mbox{(d)}$\end{tabular}}}}%
    \put(0.73793711,0.37964164){\color[rgb]{0,0,0}\makebox(0,0)[lt]{\lineheight{1.25}\smash{\begin{tabular}[t]{l}$\mbox{(f)}$\end{tabular}}}}%
    \put(0,0){\includegraphics[width=\unitlength,page=1]{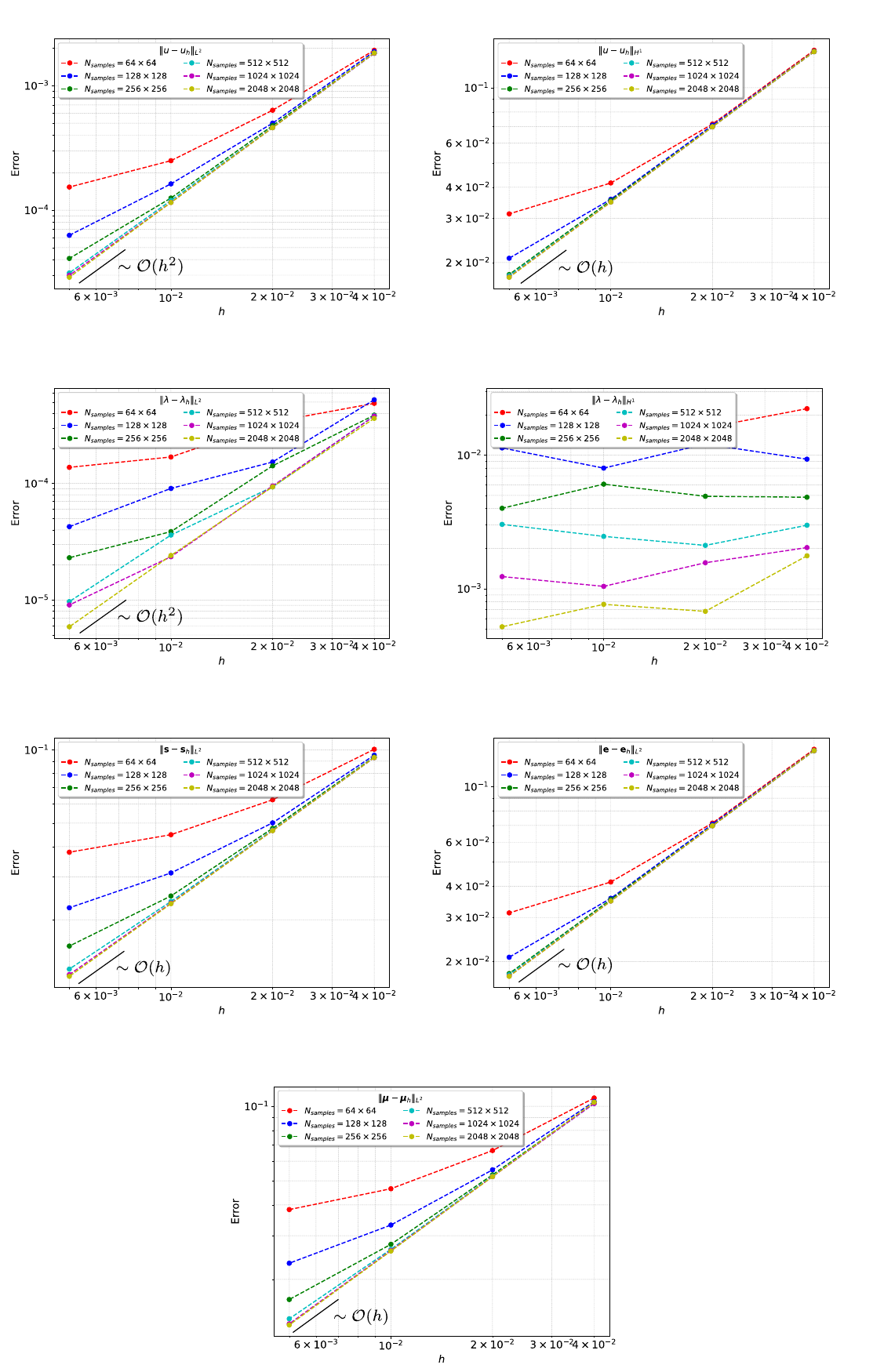}}%
    \put(0.62630892,0.06924136){\color[rgb]{0,0,0}\makebox(0,0)[lt]{\lineheight{1.25}\smash{\begin{tabular}[t]{l}$\mbox{(g)}$\end{tabular}}}}%
  \end{picture}%
\endgroup%

	\label{fig:quadricula_convergence_natural_linear}
\end{figure}
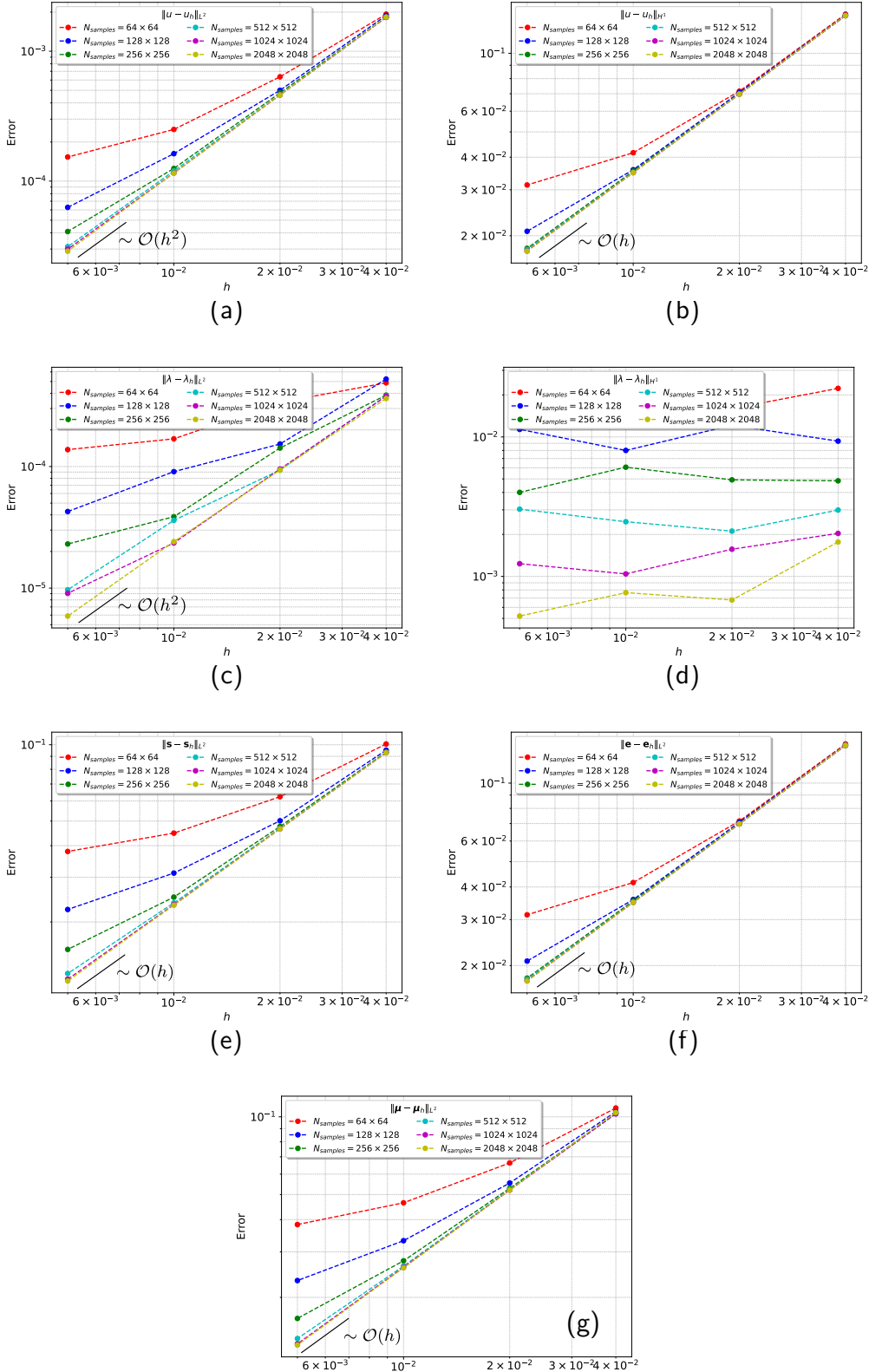

\begin{figure}[h!]
	\centering
	\caption{Mesh convergence results for the scalar fields $\uP$ and $\lP$ in the $L^2$ and $H^1$ norms, and for the vector fields $\eP$, $\sP$, and $\mP$ in the $L^2$ norm, using the  Unstabilized equal-order formulation across various numbers of subdomains $N_{\mbox{\tiny samples}}$.}
	\def\svgwidth{1.0\columnwidth}
\begingroup%
  \makeatletter%
  \providecommand\color[2][]{%
    \errmessage{(Inkscape) Color is used for the text in Inkscape, but the package 'color.sty' is not loaded}%
    \renewcommand\color[2][]{}%
  }%
  \providecommand\transparent[1]{%
    \errmessage{(Inkscape) Transparency is used (non-zero) for the text in Inkscape, but the package 'transparent.sty' is not loaded}%
    \renewcommand\transparent[1]{}%
  }%
  \providecommand\rotatebox[2]{#2}%
  \newcommand*\fsize{\dimexpr\f@size pt\relax}%
  \newcommand*\lineheight[1]{\fontsize{\fsize}{#1\fsize}\selectfont}%
  \ifx\svgwidth\undefined%
    \setlength{\unitlength}{529.13387269bp}%
    \ifx\svgscale\undefined%
      \relax%
    \else%
      \setlength{\unitlength}{\unitlength * \real{\svgscale}}%
    \fi%
  \else%
    \setlength{\unitlength}{\svgwidth}%
  \fi%
  \global\let\svgwidth\undefined%
  \global\let\svgscale\undefined%
  \makeatother%
  \begin{picture}(1,1.58289409)%
    \lineheight{1}%
    \setlength\tabcolsep{0pt}%
    \put(0.2373283,1.18811693){\color[rgb]{0,0,0}\makebox(0,0)[lt]{\lineheight{1.25}\smash{\begin{tabular}[t]{l}$\mbox{(a)}$\end{tabular}}}}%
    \put(0.73968394,1.18811693){\color[rgb]{0,0,0}\makebox(0,0)[lt]{\lineheight{1.25}\smash{\begin{tabular}[t]{l}$\mbox{(b)}$\end{tabular}}}}%
    \put(0.23794748,0.78529349){\color[rgb]{0,0,0}\makebox(0,0)[lt]{\lineheight{1.25}\smash{\begin{tabular}[t]{l}$\mbox{(c)}$\end{tabular}}}}%
    \put(0.73968394,0.78529349){\color[rgb]{0,0,0}\makebox(0,0)[lt]{\lineheight{1.25}\smash{\begin{tabular}[t]{l}$\mbox{(d)}$\end{tabular}}}}%
    \put(0.23730431,0.38247014){\color[rgb]{0,0,0}\makebox(0,0)[lt]{\lineheight{1.25}\smash{\begin{tabular}[t]{l}$\mbox{(e)}$\end{tabular}}}}%
    \put(0.74246314,0.38247014){\color[rgb]{0,0,0}\makebox(0,0)[lt]{\lineheight{1.25}\smash{\begin{tabular}[t]{l}$\mbox{(f)}$\end{tabular}}}}%
    \put(0,0){\includegraphics[width=\unitlength,page=1]{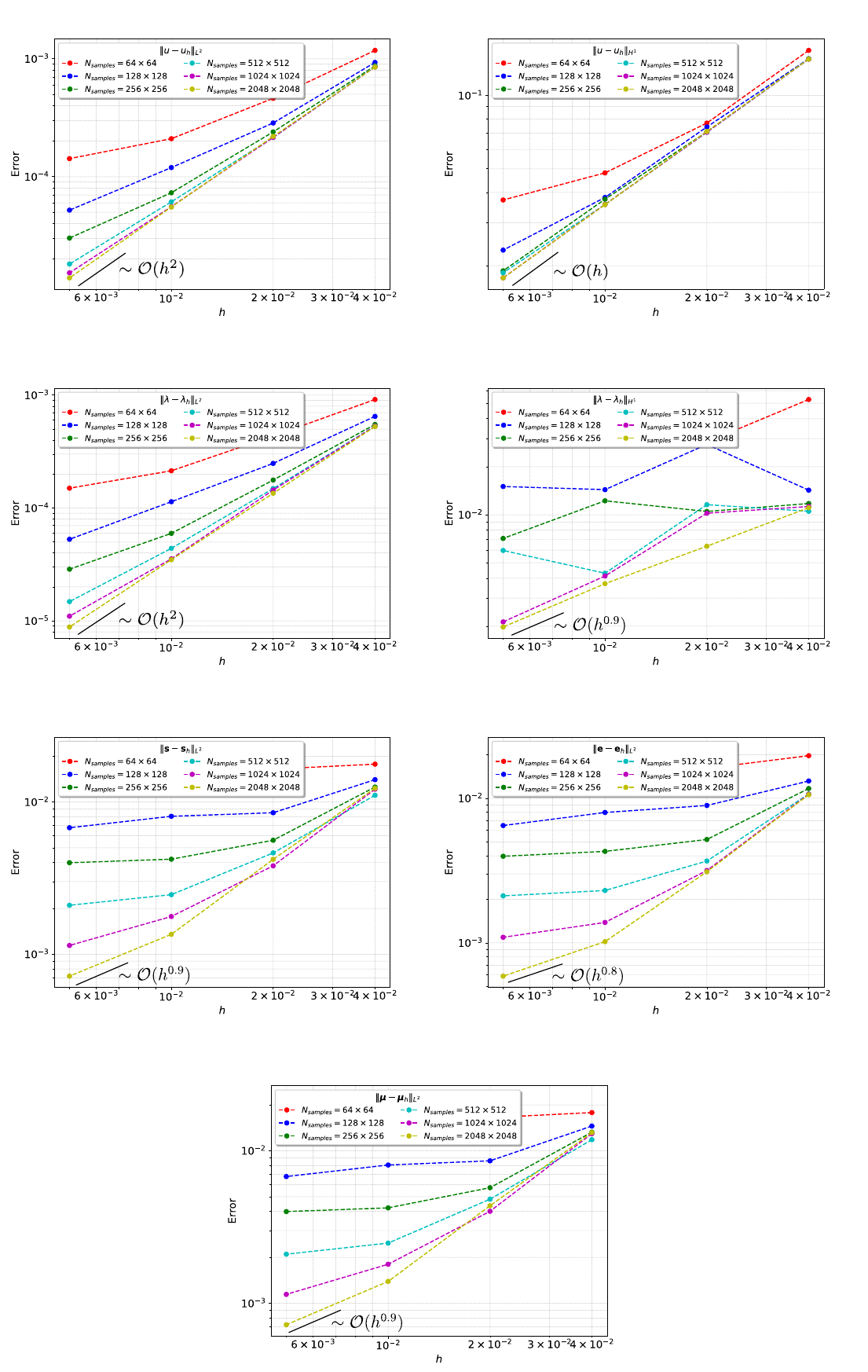}}%
    \put(0.63122175,0.07071314){\color[rgb]{0,0,0}\makebox(0,0)[lt]{\lineheight{1.25}\smash{\begin{tabular}[t]{l}$\mbox{(g)}$\end{tabular}}}}%
  \end{picture}%
\endgroup%

	\label{fig:quadricula_convergence_equal_order_unstabilized_linear}
\end{figure}

\begin{figure}[h!]
	\centering
	\caption{Mesh convergence results for the scalar fields $\uP$ and $\lP$ in the $L^2$ and $H^1$ norms, and for the vector fields $\eP$, $\sP$, and $\mP$ in the $L^2$ norm, using the Fully Stabilized equal-order formulation across various numbers of subdomains $N_{\mbox{\tiny samples}}$.}
	\def\svgwidth{1.0\columnwidth}
	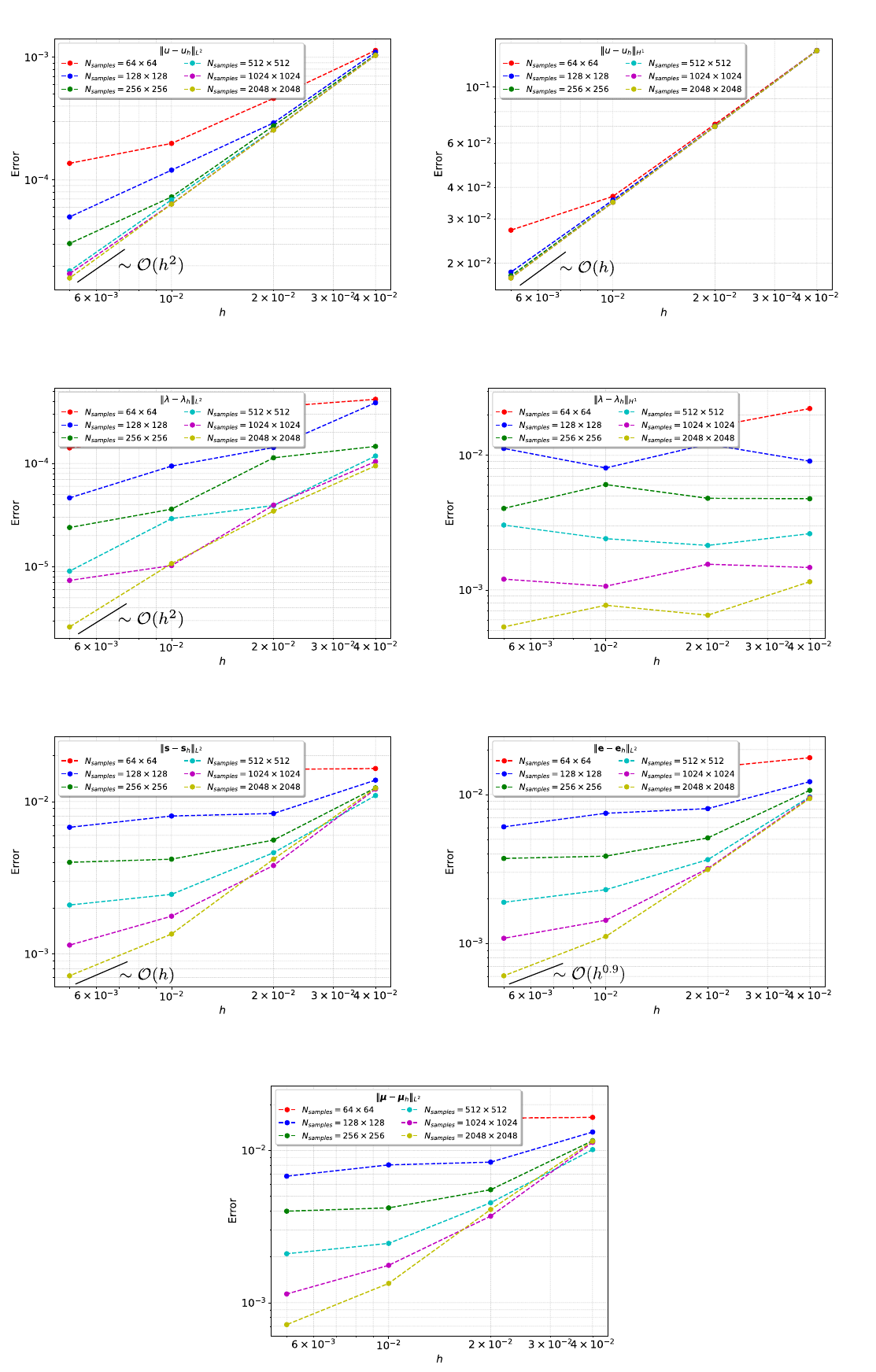
	\label{fig:quadricula_convergence_equal_order_fully_stabilized_linear}
\end{figure}

\begin{figure}[h!]
	\centering
	\caption{Mesh convergence results for the scalar fields $\uP$ and $\lP$ in the $L^2$ and $H^1$ norms, and for the vector fields $\eP$, $\sP$, and $\mP$ in the $L^2$ norm, using the Minimally Stabilized equal-order formulation across various numbers of subdomains $N_{\mbox{\tiny samples}}$.}
	\def\svgwidth{1.0\columnwidth}
	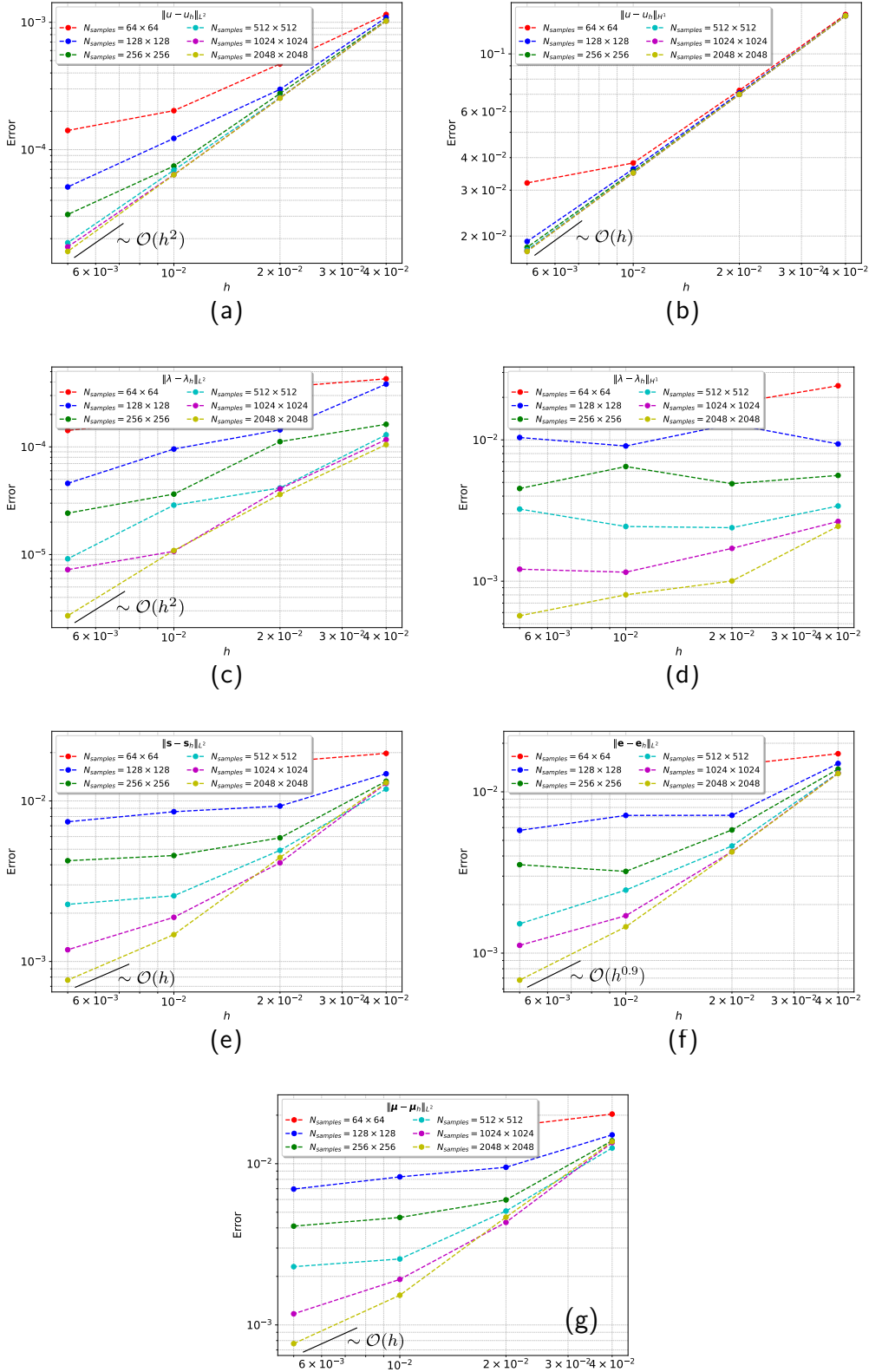
	\label{fig:quadricula_convergence_equal_order_minimally_stabilized_linear}
\end{figure}

Now, we quantitatively assess the convergence of the method under these conditions for
increasingly refined sampling partitions, from $64 \times 64$ up to $2048 \times 2048$.
The domain is discretized with a unstructured mesh of characteristic size $h$ which
is also successively refined, such that there is no conformity between the sampling
subdomains and the mesh elements. Figures~\ref{fig:quadricula_convergence_natural_linear} to~\ref{fig:quadricula_convergence_equal_order_minimally_stabilized_linear} show the mesh convergence behavior for the scalar fields $\uP$ and $\lP$ in the $L^{2}$- and $H^{1}$-norms, as well as for the vector fields $\eP$, $\sP$, and $\mP$ in the $L^{2}$-norm, for different values of $N_{\mbox{\tiny samples}}$. Figure~\ref{fig:quadricula_convergence_natural_linear} corresponds to the Natural formulation, whereas Figures~\ref{fig:quadricula_convergence_equal_order_unstabilized_linear} to~\ref{fig:quadricula_convergence_equal_order_minimally_stabilized_linear} correspond to the equal-order formulation variants: Unstabilized, Fully Stabilized, and Minimally Stabilized, respectively.

For the Natural formulation, the mesh convergence orders approach the theoretical error estimates as the number of sample subdomains $N_{\mbox{\tiny samples}}$ increases. A similar trend is observed for the variants of the equal-order formulation. In particular, the field $\uP$ exhibits $\mathcal{O}(h^{2}) $ in the $L^{2}$-norm and $ \mathcal{O}(h) $ in the $H^{1}$-norm in all formulations, with only small variations in the error magnitude.

Similarly, nearly linear convergence rates are observed for the vector fields $\sP$, $\eP$, and $\mP$ in all formulations, with the Natural formulation achieving consistent linear convergence for all three fields. Linear convergence is further achieved for the vector fields with $N_{\mbox{\tiny samples}} = 512 \times 512$ in the Natural formulation, while approximately $N_{\mbox{\tiny samples}}  = 2048 \times 2048$ are required to observe the same convergence behavior for the equal-order variants. Despite this, the equal-order formulations tend to exhibit smaller error magnitudes.

The first point to note is that the number of data samples limits the convergence rate of the method. Specifically, for a fixed number of data points, refining the finite element mesh eventually leads to error stagnation. However, as the number of samples increases, the method begins to exhibit convergence, at least within the range of $h$ values considered.

As a final observation, the field $\lP$ does not achieve the expected convergence rates in $ H^{1} $-norm, except in the Unstabilized equal-order variant. However, increasing $N_{\mbox{\tiny samples}}$ systematically reduces the associated errors, indicating an overall improvement in the approximation quality.

These results already show a consistency with respect to mesh refinement and size of the sample set for the different finite element formulations assessed which makes them a suitable candidates for data-driven applications.

\section{Conclusions}

This work presents a finite element discretization of an optimization problem arising in the emerging field of Data-Driven Computational Mechanics.
We consider the balance law associated with a diffusion-reaction problem and seek a potential field $u$ (e.g., a concentration or temperature), along with its gradient and corresponding flux $(\eP, \sP)$.
Within the Data-Driven framework, the constitutive relationship between gradients and fluxes is represented by functions $(\etP, \stP)$ derived from experimental data. As a result, these data pairs may not satisfy the governing conservation laws. The Data-Driven problem then consists of finding the fields $(u, \eP, \sP)$ that are closest to the given data, while satisfying the compatibility and equilibrium constraints.

We begin by analyzing the well-posedness of the optimization problem in the continuous setting. Subsequently, we propose a family of finite element discretizations to approximate them.
On the one hand, we consider a natural formulation in which the scalar fields are approximated using continuous piecewise polynomial functions, while the vector fields are represented by discontinuous piecewise polynomials of one degree lower.
On the other hand, we investigate equal-order formulations, both with and without stabilization, which allow the use of the same continuous interpolation for all fields.

The convergence properties of the discrete formulations were assessed through numerical experiments in 2D, using manufactured solutions to compute the error in both convex and non-convex domains. In all these experiments, the data fields $(\etP, \stP)$ were assumed to be known and available at every point in the domain. However, in a final numerical test, we also examined the consistency of the formulation with respect to the number of available data points, which better reflects scenarios typically encountered in data-driven applications.

This work represents a first step toward solving the full optimization problem that arises in Data-Driven Computational Mechanics (DDCM), where the assignment of experimental data throughout the domain must also be determined. In this context, several directions can be pursued to address the resulting mixed-integer nonlinear program (MINLP) when computed with FEM. One promising approach involves alternating direction schemes, in which data is first assigned based on the current solution, followed by solving the optimality conditions to update the fields, and repeating this process iteratively. Extensions of the method to other physical problems, such as elasticity and flow in porous media, are also the subject of ongoing research.

\section{Acknowledments}

RFA and GCB thankfully acknowledge financial support from the São Paulo Research Foundation (FAPESP) under grants 2013/07375-0 and 2023/14427-8, and from the National Council for Scientific and Technological Development - CNPq under grants 308704/2022-3 and 309786/2021-5.
PBB thankfully acknowledges the support provided by the Coordenação de Aperfeiçoamento de Pessoal de Nível Superior – Brazil (CAPES), Finance Code 88887.975792/2024-00.
CGG gratefully acknowledges the financial support from the European Research Council through the ERC Consolidator Grant ``DATA-DRIVEN OFFSHORE'' (Project ID 101083157).

\appendix
	
	\section{Firedrake Implementation}
	
	The different formulations introduced in the previous sections have 
	been implemented using the finite element platform Firedrake \cite{FiredrakeUserManual}. 
	This implementation is significantly simplified by the availability of automatic differentiation tools. The approach involves defining the discrete Lagrangian in 
	equation (\ref{eq:discLag}) and computing its variations using a \texttt{derivative} 
	function.

	For clarity, the implementation is organized into several sections. 
	Some details are omitted for simplicity, but the complete code 
	is provided in \url{https://github.com/pedrobbazon/DDCM.git}
	First, all physical and numerical parameters must be properly initialized:
\begin{lstlisting}[language=Python]
# Geometry
Nx, Ny = 20, 20
Lx, Ly = 1.0, 1.0
msh = RectangleMesh(Nx, Ny, Lx, Ly, quadrilateral=False)
x = SpatialCoordinate(msh) 
h = CellDiameter(msh)

# Parameters and Manufacuted fields  
pars = {'zeta':  fd.Constant(1.0), 
        'kappa': fd.Constant(1.0),
        'alpha': fd.Constant(0.125),   
        'gamma': fd.Constant(0.125),
        'eta':   fd.Constant(0.5), 
        'beta':  fd.Constant(0.5), 'lmu': h,
        'theta': fd.Constant(0.5), 'ls': h,
        'formulation': 'Equal-Order'
    }
\end{lstlisting}
	The keyword \texttt{formulation} can be set to 
	\texttt{\textquotesingle Natural\textquotesingle}, 
	and the associated parameters should be adjusted according 
	to the chosen formulation.
	Additionally, the user must define two functions representing 
	the fields $\etP$ and $\stP$, as well as the source term $q$ 
	in the balance equation. That is,
\begin{lstlisting}[language=Python]
def etilde(x):
    # Program your function of x[0] and x[1]
    etx = ...
    ety = ...   
    return as_vector([etx, ety])

def stilde(x):
    # Program your function of x[0] and x[1]
    stx = ...
    sty = ...
    return as_vector([stx, sty])

def q(x):
    # Program your function of x[0] and x[1]
    source = ...
    return source
\end{lstlisting}	
	This implementation can be used to solve problems similar to Test Cases 1 and 2. 
	However, to handle more general scenarios, such as in Test Case 3,
	an appropriate data distribution strategy must be implemented.
	The UFL language \cite{Alnaes2014} provides an expressive way to 
	define the operators involved in the variational formulation. 
	The first step is to specify the finite element spaces, define the 
	symbols corresponding to the variables $\sP$, $\uP$, $\eP$, $\lP$, and $\mP$, 
	and express the Lagrangian function:
\begin{lstlisting}[language=Python]
 # Lagrangian:
 if(formulation == 'Natural'):
   SF = FunctionSpace(msh, family='CG', degree=1)
   VF = VectorFunctionSpace(msh, family='DG', degree=0)
 elif(formulation == 'Equal-Order'):
   SF = FunctionSpace(msh, family='CG', degree=1)
   VF = VectorFunctionSpace(msh, family='CG', degree=1)

 # Mixed space:
 #   s    u    e    lam  mu
 W = VF * SF * VF * SF * VF

 # Unknown field
 w = Function(W)
 s, u, e, lam, mu = split(w)

 kap = pars['kappa'] # assert units in objective function
 # Proximity function
 dgrad, dflux = e-etilde(x), s-stilde(x)
 L = 0.5*inner(dgrad,dgrad)*dx + 0.5*kap*inner(dflux,dflux)*dx

 # Balance and Compatibility constraint
 L -= (inner(s, grad(lam)) + inner(pars['zeta']*u - q(x), lam))*dx
 L += inner(e - grad(u), mu)*dx
    
 # Stabilization terms
 if(pars['formulation'] == 'Equal-Order'): 
   CT = e - grad(u)
   L += 0.5*pars['alpha']*inner(CT,CT)*dx

   ET = e + mu - gt
   L -= 0.5*pars['gamma']*inner(ET,ET)*dx

   ST = kap*(s-st) - grad(lam)
   L -= (0.5*pars['eta']/kap)*inner(ST,ST)*dx

   BT = div(s) + pars['zeta']*u - q(x)
   L += (0.5*pars['theta']*kap*pars['ls']**2)*inner(BT,BT)*dx

   MT = div(mu) + pars['zeta']*lam
   L -= (0.5*pars['beta']*pars['lmu']**2)*inner(MT,MT)*dx
\end{lstlisting}
	Boundary conditions must also be defined. In this simple example, 
	homogeneous Dirichlet boundary conditions are imposed on the entire 
	boundary. However, more complex cases, such as mixed Dirichlet–Neumann 
	conditions, can be easily incorporated.
	Finally, the problem is solved using a PETSc solver provided by 
	the Firedrake platform:
\begin{lstlisting}[language=Python]
 # Test function
 dw = TestFunction(W) 
 # Optimallity conditions
 dL = derivative(L, w, dw)
 # Jacobian
 Jac = derivative(dL, w)
 # solve problem
 bclist = [DirichletBC(W.sub(1), Constant(0.0), 'on_boundary'), 
           DirichletBC(W.sub(3), Constant(0.0), 'on_boundary')]
 solve(dL == 0, w, bcs=bclist, J=Jac)
\end{lstlisting}
	In this implementation, we have omitted the artificial forcing term $f$, introduced in equation (\ref{eq:balance_with_f}), for simplicity. This term was included solely for 
	the purpose of verifying the method against manufactured solutions. Last but not least, notice that the implementation can be easily modified to account for nonlinear effects, 
	for instance, when the reaction term depends nonlinearly on $\uP$.

\bibliographystyle{plain}
\bibliography{bibliography}

\begin{thebibliography}{10}

\bibitem{Alnaes2014}
M.~S. Alnæs, A.~Logg, K.~B. Ølgaard, M.~E. Rognes, and G.~N. Wells.
\newblock Unified form language: A domain-specific language for weak
  formulations of partial differential equations.
\newblock {\em ACM Transactions on Mathematical Software}, 40:1--37, 2014.

\bibitem{Ammari2004}
H.~Ammari and H.~Kang.
\newblock {\em Reconstruction of Small Inhomogeneities from Boundary
  Measurements}.
\newblock Springer, 2004.

\bibitem{Babuska1970}
I.~Babuška.
\newblock Finite element method for domains with corners.
\newblock {\em Computing}, 6:264--273, 1970.

\bibitem{Babuska1972}
I.~Babuška and M.~B. Rosenzweig.
\newblock A finite element scheme for domains with corners.
\newblock {\em Numerische Mathematik}, 20:1--21, 1972.

\bibitem{Badia2010}
S.~Badia and R.~Codina.
\newblock Stabilized continuous and discontinuous {G}alerkin techniques for
  {D}arcy flow.
\newblock {\em Computer Methods in Applied Mechanics and Engineering},
  199:1654--1667, 2010.

\bibitem{Bochev2006}
P.~Bochev and C.~Dohrmann.
\newblock A computational study of stabilized, low-order {$C^0$} finite element
  approximations of {D}arcy equations.
\newblock {\em Computational Mechanics}, 38:323--333, 2006.

\bibitem{brezzi-fortin}
F.~Brezzi and M.~Fortin.
\newblock {\em {Mixed and Hybrid Finite Element Methods}}.
\newblock Springer, 1991.

\bibitem{Burman2023}
E.~Burman, P.~Hansbo, and M.~Larson.
\newblock The augmented {L}agrangian method as a framework for stabilised
  methods in computational mechanics.
\newblock {\em Archives of Computational Methods in Engineering},
  30:2579--2604, 2023.

\bibitem{Conti2018}
S.~Conti, S.~M{\"{u}}ller, and M.~Ortiz.
\newblock Data-driven problems in elasticity.
\newblock {\em Archive for Rational Mechanics and Analysis}, 229:79--123, 2018.

\bibitem{Conti2020}
S.~Conti, S.~M{\"{u}}ller, and M.~Ortiz.
\newblock Data-driven finite elasticity.
\newblock {\em Archive for Rational Mechanics and Analysis}, 237:1--33, 2020.

\bibitem{Demkowicz2006}
L.~Demkowicz.
\newblock Babu{\v{s}}ka $\iff$ {B}rezzi.
\newblock {\em ICES Report}, pages 06--08, 2006.

\bibitem{ern-guermond}
A.~Ern and J.L. Guermond.
\newblock {\em {Theory and Practice of Finite Elements}}.
\newblock Springer, 2004.

\bibitem{Galetzka2020}
A.~Galetzka, D.~Loukrezis, and H.~{De Gersem}.
\newblock Data-driven solvers for strongly nonlinear material response.
\newblock {\em International Journal for Numerical Methods in Engineering},
  122:1538--1562, 2021.

\bibitem{Gebhardt2024a}
C.~Gebhardt, S.~Lange, and M.~Steinbach.
\newblock Formulating and heuristic solving of contact problems in hybrid
  data-driven computational mechanics.
\newblock {\em Communications in Nonlinear Science and Numerical Simulation},
  134:108031, 2024.

\bibitem{Gebhardt2024b}
C.~Gebhardt, S.~Lange, and M.~Steinbach.
\newblock On the mathematical structure and numerical solution of
  discrete-continuous optimization problems in {DDCM}.
\newblock {\em \textit{Zeitschrift f{\"u}r Angewandte Mathematik und Physik} -
  ZAMP}, 2025.

\bibitem{Gebhardt:2020a}
C.~G. Gebhardt, D.~Schillinger, M.~C. Steinbach, and R.~Rolfes.
\newblock A framework for data-driven structural analysis in general elasticity
  based on nonlinear optimization: The static case.
\newblock {\em Computer Methods in Applied Mechanics and Engineering}, 365,
  2020.

\bibitem{Gebhardt:2020b}
C.~G. Gebhardt, M.~C. Steinbach, D.~Schillinger, and R.~Rolfes.
\newblock A framework for data-driven structural analysis in general elasticity
  based on nonlinear optimization: The dynamic case.
\newblock {\em International Journal for Numerical Methods in Engineering},
  121:5447--5468, 2020.

\bibitem{Gurtin2013}
M.~E. Gurtin, E.~Fried, and L.~Anand.
\newblock {\em The mechanics and thermodynamics of continua}.
\newblock Cambridge University Press, 2013.

\bibitem{FiredrakeUserManual}
D.~A. Ham, P.~H.~J. Kelly, L.~Mitchell, C.~J. Cotter, R.~C. Kirby, K.~Sagiyama,
  N.~Bouziani, S.~Vorderwuelbecke, T.~J. Gregory, J.~Betteridge, D.~R. Shapero,
  R.~B. Nixon-Hill, C.~J. Ward, P.~E. Farrell, P.~D. Brubeck, I.~Marsden, T.~H.
  Gibson, M.~Homolya, T.~Sun, A.~T.~T. {McRae}, F.~Luporini, A.~Gregory,
  M.~Lange, S.~W. Funke, F.~Rathgeber, G.-T. Bercea, and G.~R. Markall.
\newblock {\em Firedrake User Manual}.
\newblock Imperial College London, University of Oxford, Baylor University and
  University of Washington, 2023.

\bibitem{Kirchdoerfer2016}
T.~Kirchdoerfer and M.~Ortiz.
\newblock Data-driven computational mechanics.
\newblock {\em Computer Methods in Applied Mechanics and Engineering},
  304:81--101, 2016.

\bibitem{Kirchdoerfer2017}
T.~Kirchdoerfer and M.~Ortiz.
\newblock Data driven computing with noisy material data sets.
\newblock {\em Computer Methods in Applied Mechanics and Engineering},
  326:622--641, 2017.

\bibitem{Kirchdoerfer2018}
T.~Kirchdoerfer and M.~Ortiz.
\newblock Data-driven computational dynamics.
\newblock {\em International Journal for Numerical Methods in Engineering},
  113:1697--1710, 2018.

\bibitem{Masud2002}
A.~Masud and T.~Hughes.
\newblock A stabilized mixed finite element method for {D}arcy flow.
\newblock {\em Computer Methods in Applied Mechanics and Engineering},
  191:4341--4370, 2002.

\bibitem{Ponti2021}
L.~Ponti, S.~Perotto, and L.~Sangalli.
\newblock A {PDE}-regularized smoothing method for space–time data over
  manifolds with application to medical data.
\newblock {\em International Journal for Numerical Methods in Biomedical
  Engineering}, 38:e3650, 2022.

\bibitem{Reddy2013}
J.~N. Reddy.
\newblock {\em An Introduction to Continuum Mechanics}.
\newblock Cambridge University Press, 2013.

\end{thebibliography}

\end{document}